\documentclass [11pt]{article}
\newfont{\bbb} {msbm10}
\newcommand{\Bbb}[1]{\mbox{\bbb#1}}
\newcommand{\R}{\Bbb{R}}
\newcommand{\bS}{\Bbb{S}}

\newcommand{\Z}{\Bbb{Z}}

\newcommand{\D}{\Bbb{D}}

\newcommand{\sm}{\setminus}
\newcommand{\sbs}{\subset}
\newcommand{\ra}{\rightarrow}
\newcommand{\hH}{\hat H}
\newcommand{\bH}{\bar H}

\newcommand{\cP}{{\cal{P}}}

\newcommand{\tQ}{\tilde{Q}}
\newcommand{\ctQ}{\overline {(\tilde{Q})}}
\newcommand{\cQ}{\overline{Q}}
\newcommand{\tS}{\tilde{S}}
\newcommand{\tid}{\tilde{d}}

\newcommand{\tT}{\tilde{T}}
\newcommand{\tN}{\tilde{N}}
\newcommand{\tA}{\tilde{A}}
\newcommand{\tW}{\tilde{W}}

\newcommand{\tv}{\tilde{v}}
\newcommand{\tw}{\tilde{w}}
\newcommand{\p}{\partial}
\newcommand{\bo}{\partial_\infty}

\newcommand{\met}{{\cal{MET}}}
\newcommand{\mo}{{\cal{MET}}^{sec< 0}}
\newcommand{\mep}{{\cal{MET}}^{sec< \epsilon}}

\newcommand{\cC}{{\cal{C}}}

\usepackage[all]{xy}

\begin{document}

\title{On the Topology of the Space of Negatively Curved Metrics}
\author{F. T. Farrell and P. Ontaneda\thanks{The first author was
partially supported by a NSF grant.
The second author was supported in part 
by research  grants from CNPq(Brazil) and NSF.}}
\date{}

\maketitle

\begin{abstract} We show that the space of negatively curved metrics of a closed
negatively curved Riemannian $n$-manifold, $n\geq 10$, is highly
non-connected. 
\end{abstract}
\vspace{.3in}

\noindent {\bf \large  Section 0. Introduction.}\\

Let $M$ be a closed smooth manifold. We denote by $\met (M)$ the space of all smooth Riemannian metrics
on $M$ and we consider $\met (M)$ with the smooth topology. 
Note that the space $\met ( M ) $ is contractible. 
A subspace of metrics whose sectional curvatures lie in some interval (closed, open,
semi-open) will be denoted by placing a superscript on $\met (M)$. For example,
$\mep (M)$ denotes the subspace of $\met (M)$ of all Riemannian metrics
on $M$ that have all sectional curvatures less that $\epsilon$.
Thus saying that all sectional curvatures of a Riemannian metric $g$ lie in the interval $[a,b]$ is equivalent to
saying that $g\in\met^{\, a\,\leq\, sec\,\leq \, b}(M)$. Note that if $I\sbs J$ then
$\met^{sec\in I}(M)\sbs\met^{sec\in J}(M)$. Note also that $\met^{sec\, =\, -1}(M)$ is the space of hyperbolic
metrics ${\cal{H}}yp\, (M)$ on $M$.\\

A natural question 
about a closed negatively curved manifold $M$ is the following: is
the space $\mo (M) $ of negatively curved metrics on $M$ path connected?
This problem has been around for some time and has been posed several times in the literature.
see for instance K. Burns and A. Katok (\cite{BK}, Question 7.1).
In dimension two, Hamilton's Ricci flow \cite{H} shows that ${\cal{H}}yp\, (M^2)$
is a deformation retract of $\mo (M^2)$. But ${\cal{H}}yp\, (M^2)$
fibers over the Teichm\"{u}ller space ${\cal{T}} (M^2)\cong\R^{6\mu-6}$ ($\mu$ is the genus of $M^2$),
with contractible fiber ${\cal{D}}=\R^+\times DIFF(M^2)$ \cite{EE}.
Therefore ${\cal{H}}yp\, (M^2)$ and $\mo (M^2)$ are contractible.\\

In this paper we prove that, for $n\geq 10$,  $\mo (M^n)$ is never
path-connected; in fact, it has infinitely many path-components.
Moreover we show that all the groups $\pi_{2p-4}(\mo (M^n))$ are non-trivial
for every prime number $p>2$, and such that $p<\frac{n+5}{6}$.
(In fact, these groups contain the infinite sum $(\Z_p)^\infty$ of $\Z_p=\Z/p\,\Z$'s,
and hence they are not finitely generated).
Also, the restriction on $n=dim\, M$ can be improved to $p\leq \frac{n-2}{4}$.
See Remarks 1 below.)
We also show that $\pi_1(\mo (M^n))$ contains the infinite sum $(\Z_2)^\infty$ when $n\geq 14$.
These results about $\pi_k$ are true for each path component of $\met^{sec <0}(M^n)$;
i.e. relative to any base point.
Before we state our Main Theorem, we need some definitions.\\

Denote by  $DIFF(M)$ the group of all smooth self-diffeomorphisms of $M$.
We have that $DIFF(M)$ acts on $\met (M)$ pulling-back metrics:
$\phi g =(\phi^{-1})^*g=\phi_*g$, for $g\in\met(M)$ and $\phi\in DIFF(M)$,
that is, $\phi g $ is the metric such that $\phi :(M,g)\ra (M,\phi \, g)$ is an isometry.
Note that $DIFF(M)$ leaves invariant all spaces $\met^{sec\in I} (M)$, for any $I\sbs\R$.
For any metric $g$ on $M$ we denote by $DIFF(M)\, g$\, the orbit of $g$ by the action of 
$DIFF(M)$. We have a map $\Lambda_g: DIFF(M)\ra\met (M)$, given by $\Lambda_g(\phi )=\phi_*g$. Then
the image of $\Lambda_g$ is the orbit $DIFF(M)\, g$\, of $g$.
And $\Lambda_g$ of course naturally factors through $\met^{sec\in I}(M)$, if $g\in\met^{sec\in I}(M)$.
Note that if $dim\, M\geq 3$ and $g\in\met^{secc=-1}(M)$, then the statement of Mostow's
Rigidity Theorem is equivalent to saying that the map $\Lambda_g: DIFF(M)\ra\met^{secc=-1}(M)={\cal{H}}yp\, (M)$
is a surjection. Here is the statement of our main result.\\

\noindent {\bf Main Theorem.}
{\it Let $M$ be a closed smooth $n$-manifold and let $g$
be a negatively curved Riemannian metric on $M$.
Then we have that:}
\begin{enumerate}
\item[{i.}] {\it the map $\pi_0(\Lambda_g):\pi_0(\, DIFF(M)\, )\ra\pi_0(\, \mo (M)\, )$
is not constant, provided $n\geq 10$.}

\item[{ii.}] {\it the homomorphism $\pi_1(\Lambda_g):\pi_1(\, DIFF(M)\, )\ra\pi_1(\, \mo (M)\, )$
is non-zero, provided $n\geq 14$.}

\item[{iii.}] {\it For $k=2p-4$, $p$ prime integer and $1<k\leq\frac{n-8}{3}$, 
the homomorphism $\pi_k(\Lambda_g):\pi_k(\, DIFF(M)\, )\ra\pi_k(\, \mo (M)\, )$
is non-zero. (See Remarks 1 below.)}
\end{enumerate}

\noindent {\bf Addendum to the Main Theorem.} {\it We have that the image of
$\pi_0(\Lambda_g)$ is infinite and in cases (ii.), (iii)
mentioned in the Main Theorem, the image of $\pi_k(\Lambda_g)$ is not finitely
generated. In fact we have:}
\begin{enumerate}
\item[{i.}] {\it For $n\geq 10$, 
$\pi_0(\, DIFF(M)\, )$ contains $(\Z_2)^\infty$, and  $\pi_0(\Lambda_g)|_{(\Z_2)^\infty}$
is one-to-one.}

\item[{ii.}] {\it For $n\geq 14$, the image of $\pi_1(\Lambda_g)$ contains $(\Z_2)^\infty$.}

\item[{iii.}] {\it For $k=2p-4$, $p$ prime integer and $1<k\leq\frac{n-8}{3}$, the image of $\pi_k(\Lambda_g)$ contains $(\Z_p)^\infty$.
See Remark 1 below.}
\end{enumerate}
\vspace{.1in}

For $a<b<0$ the map $\Lambda_g$
factors through the inclusion map $\met^{\, a\, \leq\, sec\,\leq\, b}(M)\hookrightarrow\met^{sec<0}(M)$ 
provided $g\in\met^{\, a\, \leq\, sec\,\leq\, b}(M)$. Therefore we have:\\
 
\noindent {\bf Corollary 1.} {\it Let $M$ be a closed smooth $n$-manifold, $n\geq 10$.
Let $a<b<0$ and assume that $\met^{\, a\, \leq\, sec\,\leq\, b}(M)$ is not empty.
Then the inclusion map  $\met^{\, a\, \leq\, sec\,\leq\, b}(M)\hookrightarrow\met^{sec<0}(M)$ 
is not null-homotopic. Indeed, the induced maps, at the $k$-homotopy level, are not constant for $k=0$, and non-zero
for the cases (ii.), (iii.) mentioned in the Main Theorem. Furthermore, the image of these
maps satisfy a statement analogous to the one in the Addendum  to the Main Theorem.}\\

\noindent If $a=b=-1$ we have\\

\noindent {\bf Corollary 2.} {\it Let $M$ be a closed hyperbolic $n$-manifold, $n\geq 10$..
Then the inclusion map  ${\cal{H}}yp\, (M)\hookrightarrow\met^{sec<0}(M)$ 
is not null-homotopic. Indeed, the induced maps, at the $k$-homotopy level, are not constant for $k=0$, and non-zero
for the cases (ii.), (iii.) mentioned in the Main Theorem. Furthermore, the image of these
maps satisfy a statement analogous to the one in the Addendum  to the Main Theorem.}\\

Hence, taking $k=0$ (i.e. $p=2$) in Corollary 2, we get that for any closed hyperbolic  manifold $(M^n, g)$, $n\geq 10$,
there is a hyperbolic metric $g'$ on $M$ such that $g$ and $g'$ cannot be joined by a
path of negatively curved metrics. \\

Also, taking $a=-1-\epsilon$, \, $b=-1$\, ($0\leq\epsilon $)  in Corollary 1
we have that the space $\met^{ \, -1-\epsilon \,\leq \, sec\,\leq -1}(M^n)$ of $\epsilon$-pinched 
negatively curved Riemannian metrics on $M$ has infinitely many path components, provided 
it is not empty and $n\geq 10$. And
the homotopy groups $\pi_k ( \met^{\, -1-\epsilon\, \leq\, sec\,\leq -1}(M) )$, are non-zero
for the cases (ii.), (iii.) mentioned in the Main Theorem. Moreover, these groups are not finitely generated.\\
 
\noindent {\bf Remark 1.} The restriction on $n=dim\, M$ given in the Main Theorem, its Addendum and its
Corollaries are certainly not optimal. In particular, in (iii.) it can be improved to
$1<k<\frac{n-10}{2}$ by using Igusa's ``Surjective Stability Theorem'' (\cite{I}, p. 7).\vspace{.4in}

Another interesting application of the Main Theorem shows that the answer to the following
natural question is negative:\\

\noindent {\bf Question:} {\it Let $E\ra B$ be a fibre bundle whose fibres are diffeomorphic
to a closed negatively curved manifold $M^n$. Is it always possible to equip its fibres
with negatively curved Riemannian metrics (varying continuously from fibre to fibre)?}\\

The negative answer is gotten by setting $B=\bS^{k+1}$, where $k$ is as in the Main Theorem
case (iii) (or $k=0,\, 1$, case (i), (ii)), and the bundle $E\ra\bS^{k+1}$ is obtained by the
standard clutching construction using an element $\alpha\in\pi_k(DIFF(M))$ such that
$\pi_k(\Lambda_g)(\alpha)\neq 0$, for {\it every} negatively curved Riemannian metric $g$ 
on $M$. Using our method for proving the Main Theorem (in particular Theorem 1 below)
one sees that such elements $\alpha$, which are independent of $g$, exist
in all cases (i), (ii), (iii).
\vspace{.5in}

The Main Theorem follows from Theorems 1 and 2 below. Before we state these results we need some
definitions and constructions. For a manifold $N$ let $P(N)$ be the space of
topological pseuso-isotopies of $N$, that is, the space of all homeomorphisms
$N\times I\ra N\times I$, $I=[0,1]$, that are the identity on $(N\times \{ 0\})\cup (\p N\times I)$.
We consider $P(N)$ with the compact-open topology.
Also, $P^{diff}(N)$ is the space of all {\it smooth} pseudo-isotopies on $N$, with the smooth topology.
Note that $P^{diff}(N)$ is a subset of $P(N)$.
The map of spaces $P^{diff}(N)\ra P(N)$ is continuous and will be denoted by $\iota_N$, or simply by $\iota$.
The space of all self-diffeomorphisms of $N$ will be denoted by $DIFF(N)$, considered with the smooth topology. Also
$DIFF(N,\p)$ denotes the subspace of $DIFF(N)$ of all self-diffeomorphism of $N$ which are the identity on $\p N$.\\

\noindent {\bf Remark 2.} We will assume that the elements in $DIFF(N,\p)$ are the
identity {\it near} $\p N$.\\

Note that $DIFF(N\times I, \,\p\, )$
 is the subspace of $P^{diff}(N)$ of all smooth pseudo-isotopies 
whose restriction to $N\times\{ 1\} $ is the identity. The restriction of $\iota_N$ to $DIFF(N\times I,\p)$
will also be denoted by $\iota_N$.
The map $\iota_N : DIFF(N\times I,\, \p\, )\ra P(N)$ is one of the ingredients in the statement Theorem 1. \\

We will also need the following construction. Let $M$ be a negatively curved
$n$-manifold. Let $\alpha :\bS^1\ra M$ be an embedding. 
Sometimes we will denote the image $\alpha (\bS^1)$ just by $\alpha$.
We assume that the normal bundle of
$\alpha$ is orientable, hence trivial. Let $V:\bS^1\ra TM\times...\times TM$, be an
orthonormal trivialization of this bundle: $V(z)=(v_1(z),...,v_{n-1}(z))$ is
an orthonormal base of the orthogonal complement of $\alpha (z)'$ in $T_zM$.
Also, let $r>0$, such that $2r$ is less that the width of the normal geodesic tubular
neighborhood of $\alpha$. Using $V$, and the exponential map of geodesics orthogonal
to $\alpha$ we identify the normal geodesic tubular neighborhood
of width $2r$ minus $\alpha$, with $\bS^1\times \bS^{n-2}\times (0,2r]$.
Define $\Phi =\Phi^M (\alpha , V , r): DIFF(\bS^{1}\times\bS^{n-2}\times I,\p )\ra DIFF (M)$ in the following way.
For $\varphi\in DIFF(\bS^1\times\bS^{n-2}\times I,\p )$ let $\Phi (\varphi ): M\ra M$  be
the identity outside $\bS^1\times \bS^{n-2}\times [r,2r]\sbs M$, and $\Phi (\varphi)= \lambda^{-1}\varphi\lambda$,
where $\lambda(z,u,t)=(z,u,\frac{t-r}{r})$, for $(z,u,t)\in\bS^1\times \bS^{n-2}\times [r,2r]$.
Note that the dependence of $\Phi (\alpha, V, r)$ on $\alpha$ and $V$ is essential, while its
dependence on $r$ is almost irrelevant.\\

We denote by $g$ the negatively curved metric on $M$.
Hence we have the following diagram
$$\begin{array}{ccccc} DIFF( \, (\bS^1\times\bS^{n-2})\times I ,\p\, )&\stackrel{\Phi}{\ra}& DIFF(M)&
\stackrel{\Lambda_g}{\ra}&\mo (M)\\&&&&\\
\iota\,\, \downarrow \,\,\,\,\,& & &&\\&&&&\\  P(\bS^1\times\bS^{n-2})&&&&\end{array}$$

\noindent where $\iota=\iota_{ \bS^1\times\bS^{n-2}  }$ and $\Phi=\Phi^M (\alpha , V , r )$.\\

\noindent {\bf Theorem 1.} {\it Let $M$ be a closed $n$-manifold with a negatively
curved metric $g$. Let $\alpha$, $V$, $r$ and $\Phi=\Phi (\alpha , V, r)$ be as above,
and assume that $\alpha$ in not null-homotopic. Then 
$Ker \, (\, \pi_k(\Lambda_g \Phi)\,)\sbs Ker\, (\, \pi_k(\iota )\, )$, for $k<n-5$.
Here $\pi_k(\Lambda_g \Phi)$ and $ \pi_k(\iota )$ are the homomorphisms
at the $k$-homotopy group level induced by  $\Lambda_g \Phi$ and $\iota=\iota_{\bS^1\times\bS^{n-2}}$,
respectively.}\\

\noindent {\bf Remark.} In the statement of Theorem 1 above,
by $Ker \, (\, \pi_0(\Lambda_g \Phi)\,)$ (for $k=0$)
we mean the set $ \left( \pi_0(\, \Lambda_g\Phi\, )\right)^{-1}([g])$, where $[g]\in\pi_0(\met^{sec<0}(M))$ 
is the connected component of the metric $g$.\\

Hence to deduce the Main Theorem from Theorem 1 we need to know that 
$\pi_k(\iota_{\bS^1\times\bS^{n-2}})\, )$
is a non-zero homomorphism. Furthermore, to prove  the Addendum to the Main Theorem we have to
show that $\pi_k(DIFF(\bS^1\times\bS^{n-2}\times I,\p ))$ contains an infinite sum of
$\Z_p$'s (resp. $\Z_2$'s) where $k=2p-4$, $p$ prime (resp. $k=1$)
and $\pi_k(\iota_{\bS^1\times\bS^{n-2}})\, )$
restricted to this sum is one-to-one.\\

\noindent {\bf Theorem 2.} {\it Let $p$ be a prime integer such that $max\,\{ 9,6p-5\}<n$.
Then for $k=2p-4$ we have that
$\pi_k(DIFF(\bS^1\times\bS^{n-2}\times I,\p ))$ contains $(\Z_p)^\infty$ and
$\pi_k(\iota_{\bS^1\times\bS^{n-2}})\, )$ restricted to $(\Z_p)^\infty$ is one-to-one.}\\

\noindent {\bf Addendum to Theorem 2.}
{\it Assume $n\geq 14$. Then $\pi_1(DIFF(\bS^1\times\bS^{n-2}\times I,\p ))$ contains $(\Z_2)^\infty$ and
$\pi_1(\iota_{\bS^1\times\bS^{n-2}})\, )$ restricted to $(\Z_2)^\infty$ is one-to-one.}\\

The paper is structured as follows. In Section 1 we give some Lemmas, including some fibered
versions of Whitney embedding Theorem. In Section 2 we give (recall) some facts about simply connected negatively curved
manifolds and their natural extensions to a special class of non-simply connected ones. The results and facts in Sections
1 and 2 are used in the proof of Theorem 1, which is given in section 3. Finally Theorem 2 is proved in section 4.\\

Before we finish this introduction, we sketch an argument that, we hope, motivates our proof of Theorem 1. To avoid complications, let's
just consider the case $k=0$. In this situation we want to show the following:\\

\noindent {\it Let $\theta\in DIFF(\bS^1\times\bS^{n-2}\times I,\p )\sbs P(\bS^1\times\bS^{n-2})$,  and write 
$\varphi=\Phi(\theta):M\ra M$.
Suppose that $\theta$ cannot be joined by a path to the identity in $P(\bS^1\times\bS^{n-2})$. Then $g$ cannot be joined 
to $\phi_*g$ by a path of negatively curved metrics.}\\

Here is an argument that we could tentatively use to prove  the statement above. Suppose that there is a smooth path $g_u$, $u\in [0,1]$, of
negatively curved metrics on $M$, with $g_0=g$ and $g_1=\varphi_*g$. 
We will use $g_u$ to show that $\theta$ can be joined to the identity in $P(\bS^1\times\bS^{n-2})$. 
We assume that $\alpha$ is an embedded closed geodesic in $M$.
Let $Q$ be the cover of $M$ corresponding to
the infinite cyclic group generated by $\alpha$. 
Each $g_u$ lifts to a $g_u$ on $Q$ (we use the same letter).
Then $\alpha$ lifts isometrically to $(Q, g)$ and we can identify
$Q$ with $\bS^1\times\R^{n-1}$ such that $\alpha$ corresponds $\,\bS^1=\bS^1\times \{ 0\}$ and such that
each $\{ z\}\times \R v$, $v\in\bS^{n-2}\sbs\R^{n-1}$, corresponds to a $g$ geodesic ray emanating perpendicularly from
$\alpha$. For each $u$, the complete negatively curved manifold $(Q,g_u)$
contains exactly one closed geodesic $\alpha_u$, and $\alpha_u$ is freely homotopic to $\alpha$. Let us assume
that $\alpha_u=\alpha$, for all $u\in [0,1]$. Moreover, let us assume that $g_u$ coincides with $g$ in the
normal tubular neighborhood $W$ of length one of $\alpha$. Note that $Q\sm\, int\, W$ can be identified with
$(\bS^1\times\bS^{n-2})\times [1,\infty)$.
Using geodesic rays emanating perpendicularly from
$\alpha$, we can define a path of diffeomorphisms $f_u:(\bS^1\times\bS^{n-2})\times [1,\infty)\ra
(\bS^1\times\bS^{n-2})\times [1,\infty)$ by 
$f_u= [exp]^{-1}\circ exp^u $, where $exp^u$ denotes the normal (to $\alpha$) exponential
map with respect to $g_u$, and $exp=exp^0$. Using ``the space at infinity'' $\bo Q$ of $Q$ (see Section 2) we can extend
$f_u$ to $(\bS^1\times\bS^{n-2})\times [1,\infty]$ which we identify with $(\bS^1\times\bS^{n-2})\times [0,1]$.
Finally, it is proved that $f_1$ can be joined to $\theta$ in $P(\bS^1\times\bS^{n-2})$
(see Claim 6 in Section 3). This is enough because $f_0$ is the 
identity. 

Along the ``sketch of the proof" above we have of course made several unproven claims (that will be proven later); and we have also made 
a few assumptions: (1) $\alpha $ is an embedded closed geodesic,
(2) $\alpha_u=\alpha$\, for all $u$, (3) $g_u$ coincides with $g$ in a neighborhood of $g$.
Item (1) can be obtained ``after a deformation'' in $Q$. Item (2) can also be obtained after a deformation in $Q$ using the results
of Section 2. We do not know how to obtain (3) after a deformation (and this might even be impossible to do) so
we have to use some approximation methods based on Lemma 1.6 which implies that we can take a very thin normal
neighborhood $W$ of $\alpha$ such that all normal (to $\alpha$) $g_u$ geodesics rays will intersect $\p W$ transversally
in one point.\\

\noindent{\bf Acknowledgment.} We wish to thank Tom Goodwillie for communicating the $p$-torsion
 Theorem to us. This Theorem appears at the end of Section 4 and is crucial to the proof of Theorem 2 (when $k>1$).
\vspace{.6in}

\noindent {\bf \large  Section 1. Preliminaries.}\\

For smooth manifolds $A$, $B$, with $A$ compact, $C^\infty(A,B)$, $DIFF(A)$,
$Emb \,(A,B)$ denote the space of smooth maps, smooth self-diffeomorphisms and smooth embeddings of
$A$ into $B$, respectively. We consider these spaces with the smooth topology.
The $l$-disc will be denoted by $\D^l$. We choose  $u_0=(1,0,...,0)$ as the base point
of $\bS^{l}\sbs\D^{l+1}$ . For a map $f:A\times B\ra C$, we denote by $f_a$ the map given by
$f_a(b)=f(a,b)$. A map $f:\D^l\times A\ra B$ is {\it radial near} $\p$
if $f_u=f_{tu}$ for all $u\in\p\D^l=\bS^{l-1}$ and $t\in [1/2,1]$.
Note that any map $f:\D^l\times A\ra B$ is homotopic rel $\p\D^l\times A$ to a
map that is radial near $\p$.
The next Lemma is a special case of a parametrized version of Whitney's
Embedding Theorem.\\

\noindent {\bf Lemma 1.1.} {\it Let $P^m$ and $D^{k+1}$ be compact smooth manifolds and 
let $T$ be a closed smooth submanifold of $P$.  
Let $Q$ be an open subset of $\R^n$ and let
$H':D\times P\ra Q $ be a smooth map such that:
(1) $H'_u|_T:T\ra Q$ is an embedding for all $u\in D$, \, (2) $H'_u$ is an embedding for all $u\in\p D$.
Assume that that $k+2m+1<n$. 
Then  $H'$ is homotopy equivalent to a smooth map 
$\bH:D\times P\ra Q $ such that}
\begin{enumerate}
\item[{1.}] {\it $\bH_u:P\ra Q$ is an embedding, for all $u\in D$.}

\item[{2.}] {\it $\bH|_{D\times T}=H'|_{D\times T}$.}

\item[{3.}]  {\it $\bH|_{\p D\times P}=H'|_{\p D\times P}$.}
\end{enumerate}

\noindent {\bf Proof.} It is not difficult to construct a smooth map $g:P\ra\R^q$ such that:
(i) $g:P\sm T\ra\R^q\sm\{ 0\}$ is a smooth embedding \,\, (ii) $g(T)=\{ 0\}\in\R^q$\,\,
(iii) $D_p\, g\, ( v)\neq0$, for every $p\in T$ and $v\in T_pP\sm T_pT$.
Let $\varpi:D\ra [0,1]$ be a smooth map such that $\varpi^{-1}(0)=\p D$.
Define $G=H'\times g:D\times P\ra Q\times \R^q$,
$G(u,p)=(H'(u,p), \varpi (u)g(p))$. Then, for each $u\in D$, $G_u : P\ra Q\times \R^q$ is an embedding.
Moreover, $G|_{D\times T}=H'_{D\times T}$, where we consider $Q=Q\times\{ 0\}\sbs Q\times \R^q$.
Also $G|_{\p D\times P}=H'|_{\p D\times P}$.
Note that $G$ is homotopic to $H'$ because $g$ is homotopy trivial.
Now, as in the proof of Whitney's Theorem, we want to reduce
the dimension $q$ to $q-1$. So assume $q>0$.
 Given $w\in \bS^{n+q-1}\sbs\R^{n+q}=\R^n\times\R^q$,\, 
$w\notin\R^n\times\R^{q-1}=\R^{n+q-1}$, denote by $L_w:\R^{n+q}\ra\R^{n+q-1}$ the linear projection
``in the $w$-direction".
As in the proof of Whitney's Theorem, using the dimension restriction and Sard's Theorem, we
can find a ``good" $w$:\\

\noindent {\bf Claim.} {\it There is a $w$ such that $L_w|_{G_u(P)}:G_u(P)\ra\R^{n+q-1}$ is an embedding, for all 
$u\in D$.}\\

\noindent For this consider:
$$\begin{array}{ll}
r:D\times (\, (P\times P)\sm \Delta(P)\, )\ra\R^{n+q}, & \,\, r(u,p,q)=\frac{H'_u(p)-H'_u(q)}{|H'_u(p)-H'_u(q)|}\\ \\
s:D\times SP\ra\R^{n+q}, & \,\,  s(u,v)=\frac{D_p(H'_u) (  v) }{| D_p(H'_u)( v)  |},\,\,\, v\in T_pP
\end{array}$$

\noindent Here $\Delta (P)=\{(p,p)\, :\, p\in P\}$ and $SP$ is the sphere bundle of $P$
(with respect to any metric).
Since $(k+1)+2m<n$ and $q>0$, by Sard's Theorem the images of $r$ and $s$ have measure zero
in $\bS^{n+q-1}$. This proves the Claim.

Also, since $D$ and $P$ are compact, we can choose $w$ close enough to
(0,...,0,1) such that $L_w(\, G(D\times P)\, )\sbs Q\times\R^{q-1}$. Define $G_1=L_wG$.
In the same way we define  $G_2:D\times P\ra Q\times \R^{q-2}$ and so on.
Our desired map $\bH$ is $\bH=G_q$. This proves the Lemma.\\

In what follows of this section we consider $Q=\bS^1\times\R^{n-1}=(\bS^1\times \R)\times \R^{n-2}\sbs\R^2\times\R^{n-2}$,
where the inclusion $\bS^1\times\R\hookrightarrow\R^2$ is given by $(z,s)\mapsto e^sz$. That is,
we identify $\bS^1\times\R$ with the open set $\R^2\sm\{0\}$, hence we identify $Q=\bS^1\times\R^{n-1}$
with $(\R^2\sm\{ 0\} )\times \R^{n-2}=\R^n\sm\, (\, \{ 0\}\times\R^{n-2}\, )$.
Also, identify $\bS^1$ with $\bS^1\times\{ 0\}\sbs Q$ and denote by $h_0:\bS^1\ra\bS^1\times\R^{n-1}=Q$
the inclusion.
For $t>0$ denote by $\kappa_t:\R^2\times\R^{n-2}\ra\R^2\times\R^{n-2}$ given by
$\kappa_t (a,b)=(ta,b)$. Note that $\kappa_t$ restricts to $Q=
(\R^2\sm\{ 0\} )\times \R^{n-2}$.\\

\noindent {\bf Lemma 1.2.}  {\it Let 
$h,h':\D^{k+1}\times\bS^1\ra Q $ be continuous maps such that 
$h_u$, $ h'_{u}$ are homotopy equivalent, for all $u\in\bS^{k}$.
That is there is $H':\bS^k\times\bS^1\times I\ra Q$
such that $H'(u,z,0)=h(u,z)$, $H'(u,z,1)=h'(u,z)$, for all $(u,z)\in\bS^k\times\bS^1$.
For $k=0$ also assume that the loop $h(t,1)*H'(1,1,t)*[h'(t,1)]^{-1}*[H'(-1,1,t)]^{-1}$ is null-homotopic.
Then $H'$ extends to $H':\D^{k+1}\times \bS^1\times I\ra Q$
such that $H'_u$ is a homotopy from $h_u$ to $h'_{u}$, that is}
\begin{enumerate}
\item[{1.}] {\it $H'_u|_{\bS^1\times\{ 0\}}=h_u$, for $u\in\D^{k+1}$.}

\item[{2.}] {\it $H'_u|_{\bS^1\times\{ 1\}}=h'_u$, for $u\in\D^{k+1}$.}
\end{enumerate}

\noindent {\bf Proof.} First define $H'=h$ on $\D^{k+1}\times\bS^1\times\{ 0\}$
and $H'=h'$ on $\D^{k+1}\times\bS^1\times\{ 1\}$.
Note that $H'$ is defined on $\p \, ( \D^{k+1}\times \{ 1\}\times I )$. 
Since $Q$ is aspherical, we can extend $H'$ to $\D^{k+1}\times\{ 1\}\times I$
(for $k=0$ use the assumption given in the statement of the Lemma).
$H'$ is now defined on $A=\bS^{k}\times\bS^1\times \{ 0,1\}\,\cup\, \D^{k+1}\times\{ 1\}\times I$.
Since $\D^{k+1}\times\bS^1\times I$ is obtained from $A$
by attaching a $(k+3)$-cell and $Q$ is aspherical, we can extend $H'$ to $\D^{k+1}\times\bS^1\times I $. This
proves the Lemma.\\

\noindent {\bf Lemma 1.3.}  {\it Let
$h:\D^{k+1}\times\bS^1\ra Q $ be a smooth map which is radial near $\p$. 
Assume that $h_u\in Emb\, (\bS^1, Q)$
for all $u\in\D^{k+1}$ and $h_u=h_0$, for all $u\in\bS^k$. 
For $k=0$ assume that the loop $h(u,1)$ is homotopy trivial. If $\, k+5<n$
then there is a smooth map $\hH:\D^{k+1}\times \bS^1\times I\ra Q$
such that}
\begin{enumerate}
\item[{1.}] {\it $\hH_u|_{\bS^1\times\{ 0\}}=h_u$, for $u\in\D^{k+1}$.}

\item[{2.}] {\it $\hH_u|_{\bS^1\times\{ 1\}}=h_0$, for $u\in\D^{k+1}$.}

\item[{3.}] {\it  $\hH_u$ is a smooth isotopy from $h_u$ to $h_{0}$.}

\item[{4.}] {\it $(\hH_u)_t=h_0$, for all $u\in\bS^k$ and $t\in I$. Here $(\hH_u)_t(z)=\hH(u,z,t)$.}
\end{enumerate}

\noindent {\bf Proof.} During this proof some isotopies and functions 
have to be smoothed near endpoints and boundaries. We do not do this to avoid
unnecessary technicalities.

Let $D=\D^{k+1}_{1/2}$ be the closed $(k+1)$-disc of radius 1/2.
Since $h(\D^{k+1}\times \bS^1)\sbs Q=\R^n\sm\, (\, \{ 0\}\times\R^{n-2}\, )$, we have that
$h(\D^{k+1}\times \bS^1)$ does not intersect   $ \{ 0\}\times\R^{n-2}$. Therefore the distance  $d$ from 
$h(\D^{k+1}\times \bS^1)$ to $\{ 0\}\times\R^{n-2}$ is
positive. Let $c<1$ be such that $c<d$.\\

{\bf Definition of $(\hH_u)_t$ for $t\in [1/2,1]$.}
In this case define for $u\in\bS^k$, 
$(\hH_{su})_t=\kappa_{\lambda}h_0$, 
where {\bf (1)} $\lambda=1-4(1-t)(1-s)+4(1-t)(1-s)c$ if $s\in [1/2,1]$
and {\bf (2)} $\lambda=(2t-1)+(2-2t)c$ $s\in [0,1/2]$.\\

{\bf Definition of $(\hH_su)_t$ for $t\in [0,1/2]$ and $s\in [1/2,1]$.}
Define for $u\in\bS^k$, $s\in [1/2,1]$:
$(\hH_{su})_t=\kappa_{\lambda}$, 
where $\lambda=1-4t(1-s)+4t(1-s)c$, for $t\in [0,1/2]$.\\

{\bf Definition of $(\hH_su)_t$ for $t\in [0,1/2]$ and $s\in [0,1/2]$.}
Note that $D=\{ su\,:\, u\in\bS^k,\, s\in [0,1/2]\}$.
We now want to define $\hH$ on $D\times \bS^1\times I$.
To do this first apply Lemma 1.2 taking: $h'_u=\kappa_c h_0$ for all $u\in D$, 
$H'(u,z,t)=\hH(u,z, t/2)$, for $(u,z, t)\in\bS^k\times\bS^1\times I$. 
Hence $H'$ extends to $D\times\bS^1\times I$.
Now apply Lemma 1.1 taking: $P=\bS^1\times I$, $T=\bS^1\times\{ 0,1\}$.
To apply this Lemma note that $H'_u|_{\bS^1\times\{ 0,1\}}$ is an embedding, for all $u\in D$,
because $H'_u|_{\bS^1\times\{ 0\}}=h_u$,  $H'_u|_{\bS^1\times\{ 1\}}=\kappa_c h_0$ are embeddings
and the images of $h_u$ and $\kappa_c h_0$ are disjoint (by the choice of $c$).
Let then $\bH$ be the map given by Lemma 1.1. Finally define $\hH(u,z,t)=\bH(u,z,2t)$.
This proves the Lemma.\\

Extending the isotopies $\hH_u$ between $h_u$ and $h'_u$ given in the Lemma above, 
to compactly supported ambient isotopies we obtain the following Corollary:\\

\noindent {\bf Lemma 1.4.} {\it Let
$h:\D^{k+1}\times\bS^1\ra Q $, be a smooth map which is radial near $\p$. 
Assume $h_u\in Emb(\bS^1,Q)$ for all $u\in\D^{k+1}$ and
that $h_u=h_{0}\in Emb\, (\bS^1, Q)$
for all $u\in\bS^k$, and $k+5<n$. Identify $\bS^1$ with $\bS^1\times\{ 0\}\sbs Q$.
For $k=0$ assume that the loop $h(u,1)$ is null-homotopic.
Then there is a smooth map $H:\D^{k+1}\times Q\times I\ra Q$
such that}
\begin{enumerate}
\item[{1.}] {\it $H_u|_{\bS^1\times\{ 0\}}=h_u$, for $u\in\D^{k+1}$.}

\item[{2.}] {\it $H_u|_{\bS^1\times\{ 1\}}=h_{0}$, for $u\in\D^{k+1}$.}

\item[{3.}] {\it  $H_u$ is a ambient isotopy from $h_u$ to $h_{0}$, that is $(H_u)_t:Q\ra Q$ 
is a diffeomorphism for all $u\in\D^{k+1}$, $t\in I$ and $(H_u)_1=1_Q$. 
Also, $H_u$ is supported on a compact subset $K\sbs Q$, where $K$ is independent of $u\in\D^{k+1}$.}

\item[{4.}] {\it $(H_u)_t=1_Q$, for all $u\in\bS^k$ and $t\in I$.}
\end{enumerate}
\vspace{.4in}

\noindent {\bf Remark.}  Note that Lemma 1.4 can be paraphrased as follows: {\it
the homotopy fiber of $Emb(\bS^1, Q)\ra C^\infty (\bS^1, Q)$ is $(n-5)$-connected.}\\

We will also need the result stated in Lemma 1.6, below. First
we prove a simplified version of it.
The $k$-sphere of radius $\delta$, $\{ v\in\R^{k+1}\, : \, |v|=\delta\}$,
will be denoted by $\bS^k(\delta )$.\\

\noindent {\bf Lemma 1.5.} {\it Let $X$ be a compact space and $f: X\ra DIFF(\R^l)$
be continuous and write $f_x:\R^l\ra\R^l$ for the image of $x$ in $DIFF(\R^l)$. 
Assume $f_x(0)=0\in\R^l$, for all $x\in X$. Then there is a $\delta_0 >0$
such that, for every $x\in X$ and $\delta\leq \delta_0$,  the map $\bS^{l-1}(\delta )\ra \bS^{l-1} $ given by 
$v\mapsto \frac{f_x(v)}{|f_x(v)|}$ is a diffeomorphism.
Moreover, the map $X\ra DIFF(\bS^{l-1}(\delta ),\, \bS^{l-1})$, given by
$x\mapsto\, (v\mapsto \frac{f_x(v)}{|f_x(v)|})$, is continuous.}\\

\noindent {\bf Proof.} First note that for all $x\in X$ and $\delta >0$, the maps in $DIFF(\bS^{l-1}
(\delta ),\, \bS^{l-1})$ given by $(v\mapsto \frac{f_x(v)}{|f_x(v)|})$ all have degree 1 or -1.
For $v\in\R^l\sm\{ 0\}$, denote by $L_x(v)$ the image of
the tangent space $T_v(\bS^{l-1}(|v|))$ by the derivative of $f_x:\R^l\ra\R^l$.
It is enough to prove that there is $\delta_0 >0$ such that
$f_x(v)\notin L_x(v)$, for all $x\in X$ and $v\in\R^l$  satisfying $0<|v|\leq \delta_0$
(because then the maps $(v\mapsto \frac{f_x(v)}{|f_x(v)|})$ would be immersions of degree 1 (or -1),
hence diffeomorphisms).\\


Suppose this does not happen. Then there is a sequence of points $(x_m,v_m)\in X\times \R^l\sm\{ 0\}$ with
\begin{enumerate}
\item[{a.}] $v_m\ra 0$.

\item[{b.}] $f_{x_m}(v_m)\in L_{x_m}(v_m)$.
\end{enumerate}

\noindent Write $w_m=\frac{v_m}{|v_m|}\in\bS^{l-1}$, $r_m=|v_m|$,
$f_m=f_{x_m}$ and $D_m=D_{v_m}f_m$.
We can assume that $x_m\ra x\in X$, and that $w_m\ra w\in\bS^{l-1}$.
It follows that there is an $u_m\in T_{v_m}(\bS^{l-1}(r_m))$,  $|u_m|=1$, such that
$ D_m.\, u_m$ is parallel to $f_m(v_m)$. Note that $\langle  u_m , v_m \rangle =0 $ and $D_m ( u_m)\neq 0$.
By changing the sign of $u_m$ we can assume that
$ \frac{D_m( u_m)}{|D_m( u_m)  |}=\frac{f_m(v_m)}{| f_m(v_m)  |}$. Also, we can suppose
that $u_m\ra u\in\bS^{l-1}$.\\

\noindent {\bf Claim.} {\it We have that
$\frac{f_m(v_m)}{|f_m(v_m)|}\ra\frac{D_0f_x( w)}{|D_0f_x( w)|}$, as \,\, $m\ra\infty$.}
\vspace{.1in}

\noindent {\bf Proof of the Claim.}
Since $f$ is continuous, all second order partial derivatives of the coordinate
functions of the $f_x$ at $v$, with, say, $|v|\leq1$, are bounded by some constant.
Hence there is a constant $C>0$ such that $|f_m(v_m)-D_0f_m( v_m)|=|f_m(v_m)-f_m(0)-D_0f_m( v_m)| \leq C\, |v_m|^2$,
for sufficiently large $m$. It follows that $\frac{f_m(v_m)}{|v_m|}\ra lim_{m\ra\infty}\frac{D_0f_m( v_m)}{|v_m|}
=D_0f_x( w)\neq 0$. This implies that $\frac{|f_m(v_m)|}{|v_m|}\ra |D_0f_x( w)|\neq 0$, thus
$\frac{|v_m|}{|f_m(v_m)|}\ra \frac{1}{|D_0f_x( w)|}$. Therefore
$lim_{m\ra\infty}\frac{f_m(v_m)}{|f_m(v_m)|)}=lim_{m\ra\infty}\frac{f_m(v_m)}{|v_m|}\frac{|v_m|}{|f_m(v_m)|}
=D_0f_x( w)\, \frac{1}{|D_0f_x( w)|}$. This proves the Claim.\\

But $ \frac{D_m(u_m)}{|D_m( u_m ) |}\ra\frac{D_0f_x( u)}{|D_0f_x( u)|}$, therefore
$ \frac{D_0f_x( u)}{|D_0f_x( u)|}=\frac{D_0f_x( w)}{|D_0f_x( w)|}$. This is a contradiction
since $D_0f_x$ is an isomorphism and $u,w\in\bS^{l-1}$ are linearly independent
(because $\langle u , w\rangle = lim_m \langle  u_m , \frac{v_m}{|v_m|} \rangle =0$).
This proves the Lemma.\\

\noindent {\bf Lemma 1.6.} {\it Let $X$ be a compact space, $N$ a closed
smooth manifold and $f: X\ra DIFF(N\times\R^l)$
be continuous and write $f_x=(f^1_x,f^2_x):N\times\R^l\ra N\times\R^l$ for the image of $x$ in $DIFF(N\times\R^l)$. 
Assume $f_x(z,0)=(z,0)$, for all $x\in X$ and $z\in N$, that is, $f_x|_N=1_N$, where we identify
$N$ with $N\times\{ 0\}$. Then there is a $\delta_0 >0$
such that, for every $x\in X$,  the map $N\times\bS^{l-1}(\delta )\ra N\times\bS^{l-1} $ given by 
$(z,v)\mapsto (f^1_x(z,v),\frac{f^2_x(z,v)}{|f^2_x(z,v)|})$ is a diffeomorphism for all $\delta\leq \delta_0$.
Moreover, the map $X\ra DIFF(N\times\bS^{l-1}(\delta ),\,  N\times\bS^{l-1})$, given by
$x\mapsto\, (\,\, (z,v)\mapsto (f^1_x(z,v), \frac{f^2_x(z,v)}{|f^2_x(z,v)|})\,\, )$, is continuous.}\\

\noindent {\bf Proof.} The proof is similar to the proof of the Lemma above. Here are the details.
Let $d=dim \, N$ and consider $N$ with some Riemannian metric. 
For $(z,v)\in N\times\R^l\sm\{ 0\}$, denote by $L_x(z,v)$ the image of
the tangent space $T_{(z,v)}(N\times\bS^{l-1}(|v|))$ by the derivative of $f_x$. As before
it is enough to prove that there is $\delta_0 >0$ such that
$(0,f^2_x(z,v))\notin L_x(z,v)\sbs (T_zN)\times \R^l=T_{(z,v)}(N\times\R^l)$, for all $x\in X$ and 
$(z,v)\in N\times\R^l$  satisfying \, $0<|v|\leq \delta_0$.\\
Before we prove this we have a Claim.\\

\noindent {\bf Claim 1.} {\it We have:}
\begin{enumerate}
\item[{1.}] {\it $D_{(z,0)}f_x^1 (y,0)=y$,\,  for all $z\in N$ and $y\in T_zN$.}

\item[{2.}] {\it $D_{(z,0)}f_x^2(y,u)=0$ implies that $u=0$.}
\end{enumerate}

\noindent {\bf Proof of Claim 1.} Since $f_x|_{N}=1_N$ we have that $D_{(z,0)}f_x (y,0)=(y,0)$, for all $y\in T_zN$.
Hence (1) holds.  Suppose $D_{(z,0)}f_x^2(y,u)=0$.
Write $y'=D_{(z,0)}f^1_x(y,u)$. Then 
$D_{(z,u)}f_x(y,u)=(y',0)=D_{(z,0)}f_x(y',0)$. But $D_{(z,0)}f_x$ is an isomorphism therefore 
$(y,u)=(y',0)$. This proves the Claim.\\

Suppose now that (2) does not happen. Then there is a sequence of points $(x_m,z_m,v_m)\in X\times N\times\R^l\sm\{ 0\}$ with
\begin{enumerate}
\item[{a.}] $v_m\ra 0$.

\item[{b.}] $(0,f_{x_m}^2(z_m,v_m))\in L_{x_m}(z_m,v_m)$.
\end{enumerate}

\noindent Write $w_m=\frac{v_m}{|v_m|}\in\bS^{l-1}$, $r_m=|v_m|$,
$f_m=f_{x_m}$ and $D_m^i=D_{v_m}f_m^i$, $i=1,2$.
We can assume that $x_m\ra x\in X$, $z_m\ra z$ and $w_m\ra w\in\bS^{l-1}$.
It follows that there is a $(s_m,u_m)\in T_{(z_m,v_m)}(N\times\bS^{l-1}(r_m))$,  $|s_m|^2+|u_m|^2=1$, such that:
(i) $D_m^1 (s_m, u_m )=0$\,\, (ii)
$ D_m^2 (s_m, u_m)$ is parallel to $f_m^2(z_m, v_m)$. We have that $\langle  u_m , v_m \rangle =0$. Since
$D_m=D_{v_m}f_m$ is an isomorphism, by (i),  
$D_m^2 (s_m,u_m)\neq 0$.
By changing the sign of $(s_m,u_m)$ we can assume that
$ \frac{D_m^2 (s_m,u_m)}{|D_m^2 (s_m,u_m)  |}=\frac{f_m^2(z_m,v_m)}{| f_m^2(z_m,v_m)  |}$. Also, we can suppose
that $u_m\ra u\in\R^{l}$ and $s_m\ra s\in T_zN$.\\

\noindent {\bf Claim 2.} {\it We have that
$\frac{f_m^2(z_m,v_m)}{|f_m^2(z_m,v_m)|}\ra\frac{D_{(z,0)}f_x^2 (0, w)}{|D_{(z,0)}f_x^2 (0,w)|}$, as \,\, $m\ra\infty$.}
\vspace{.1in}

\noindent {\bf Proof of the Claim.}
Since $f^2$ is continuous, all second order partial derivatives of the coordinate
functions of the $f_x^2$ at $v$, with, say, $|v|\leq1$, are bounded by some constant.
Hence there is a constant $C>0$ such that $|f_m^2(z_m,v_m)-D_{(z_m,0)}f_m^2(0,v_m)|=
|f_m^2(z_m,v_m)-f_m^2(z_m,0)-D_{(z_m,0)}f_m^2 (0,v_m)| \leq C\, |(0,v_m)|^2=|v_m|^2$,
for sufficiently large $m$. It follows that $\frac{f_m^2(z_m,v_m)}{|(0,v_m)|}\ra lim_{m\ra\infty}
\frac{D_{(z_m,0)}f_m^2 (0,v_m)}{|(0,v_m)|}
=D_{(z,0)}f_x^2 (0,w)$. Note that, by claim 1 and $w\neq 0$,  $D_{(z,0)}f_x^2 (0,w)\neq 0$.
This implies that $\frac{|f_m^2(z_m,v_m)|}{|(0,v_m)|}\ra |D_{(z,0)}f_x^2(0,w)|\neq 0$, thus
$\frac{|(0,v_m)|}{|f_m^2(z_m,v_m)|}\ra \frac{1}{|D_{(z,0)}f_x^2 (0,w)|}$. Therefore
$lim_{m\ra\infty}\frac{f_m^2(z_m,v_m)}{|f_m^2(z_m,v_m)|}=lim_{m\ra\infty}\frac{f_m^2(z_m,v_m)}{|(0,v_m)|}
\frac{|(0,v_m)|}{|f_m^2(z_m,v_m)|}
=D_{(z,0)}f_x^2 (0,w)\, \frac{1}{|D_{(z,0)}f_x^2 (0,w)|}$. This proves the Claim.\\

But $ \frac{D_m^2(s_m,u_m)}{|D_m^2 (s_m,u_m)  |}\ra\frac{D_{(z,0)}f_x^2(s,u)}{|D_{(z,0)}f_x^2(s,u)|}$, therefore
$ \frac{D_{(z,0)}f_x^2 (s,u)}{|D_{(z,0)}f_x^2(s,u)|}=\frac{D_{(z,0)}f_x^2 (0,w)}{|D_{(z,0)}f_x^2 (0,w)|}$. 
Consequently $D_{(z,0)}f_x^2(s,u)=D_{(z,0)}f_x^2 (0,w') $, where $w'=\lambda w$, for some $\lambda >0$.
Hence $D_{(z,0)}f_x^2 (s,u-w')=0$, and by Claim 1, $u=w'=\lambda w$ a contradiction because $|w|=1$
and $\langle u,w\rangle =0$. This proves the Lemma.\\
\vspace{.6in}

\noindent {\bf \large  Section 2. Space at infinity of some complete negatively curved manifolds.}\\

Let $(X_1,d_1)$ and $(X_2,d_2)$ be two metric spaces. A map $f:X_1\ra X_2$ is a quasi-isometric
embedding if there are $\epsilon\geq 0$ and $\lambda\geq 1$ such that
$\frac{1}{\lambda}\, d_1(x,y)-\epsilon\leq d_2(f(x),f(y))\leq \lambda\, d_1(x,y)+\epsilon$, for all
$x,y\in X_1$. A quasi-isometric embedding $f$ is called a quasi-isometry if
there is a constant $K\geq 0$ such that every point in $X_2$ lies in the $K$-neighborhood of the image of $f$.
A quasi-geodesic in a metric space $(X,d\, )$ is a quasi-isometric embedding $\beta:I\ra X$,
where the interval $I\sbs\R$ is considered with the canonical metric $d_{\R}(t,s)=|t-s|$. If $I=[a,\infty )$, $\beta$
is called a quasi-geodesic ray. 
If we want to specify the constants $\lambda$ and $\epsilon$ in the definitions above we will use the
prefix $(\lambda ,\epsilon)$.
It is a simple exercise to prove that the composition of a $(\lambda,\epsilon)$-quasi-isomeric embedding  
with a $(\lambda',\epsilon')$-quasi-isomeric embedding is a 
$(\lambda\lambda',\lambda'\epsilon+\epsilon')$-quasi-isomeric embedding. 
Also, if $f:X_1\ra X_2$ is a quasi-isometry and the Hausdorff distance between some
subsets $A,B\sbs X_1$ is finite, then the Hausdorff distance between $f(A)$ and $f(B)$ is also finite. 
In this paper a {\it unit speed geodesic} will always mean an isometric embedding 
with domain some interval $I\sbs\R$. Also a {\it geodesic} will mean a function $t\mapsto \alpha (\rho t)$,
where $\alpha$ is a unit speed geodesic and $\rho>0$.
Then every geodesic is a quasi-geodesic
with $\epsilon =0$, that is, a $(\lambda,0)$-quasi-geodesic, for some $\lambda$. \\

\noindent {\bf Lemma 2.1.} {\it Let $g,\, g'$ be two complete Riemannian metrics on the manifold $Q$.
Suppose there are constants $a,b>0$ such that
$a^2\leq g'(v,v)\leq b^2$ for every $v\in TQ$ with $g(v,v)=1$. Then the identity $(Q,g)\ra (Q,g')$ is a
$(\lambda,0)$-quasi-isometry, where $\lambda=max\{ \frac{1}{a},\, b\}$.}\vspace{.1in}

\noindent {\bf Proof.} The condition above implies: $a^2\, g(v,v)\leq g'(v,v)\leq b^2\, g(v,v)$,
which in turn implies $\frac{1}{b^2}\, g'(v,v)\leq g(v,v)\leq \frac{1}{a^2}\, g'(v,v)$, for all $v\in Q$.
Let $d,\, d'$ be the intrinsic metrics on $Q$ defined by $g$, $g'$, respectively.
Let $x,y\in Q$ and $\beta:[0,1]\ra Q$ be a path whose endpoints are $x,y$ and such that $d(x,y)=length_g (\beta)
=\int_0^1\sqrt{g(\beta'(t),\beta'(t))}\, dt$. Then $d'(x,y)\leq length_{g'}(\beta)=\int_0^1
\sqrt{g'(\beta'(t),\beta'(t))}\, dt\leq b\int_0^1\sqrt{g(\beta'(t),\beta'(t))}\, dt=b\, d(x,y)$.
In the same way we prove $d\leq\frac{1}{a}\, d'$. Then the identity $1_Q$ is a quasi-isometry with
$\epsilon=0$ and $\lambda=max\{ \frac{1}{a},\, b\}$.
This proves the Lemma.\\

In what remains of this section $(Q,g)$ will denote a complete Riemannian manifold 
with sectional curvatures in the interval $[c_1,c_2]$, $c_1<c_2<0$, and $S\sbs Q$
a closed totally geodesic submanifold of $Q$, such that the map $\pi_1 (S)\ra\pi_1(Q)$ is an isomorphism.
Write $\Gamma=\pi_1(S)=\pi_1(Q)$. Also, $d$ will denote the intrinsic metric on $Q$ induced by $g$.
Note that $S$ is convex in $Q$, hence $d|_S$ is also the intrinsic metric on $S$ induced by $g|_S$.
We can assume that the universal cover $\tilde S$ of $S$ is contained in the universal cover $\tilde Q$ of $Q$.
We will consider $\tQ$ with the lifted metric $\tilde g$ and the induced distance will be denoted by $\tilde d$. 
The group $\Gamma$ acts by isometries on $\tilde Q$ such that $\Gamma (S)=S$ and $Q={\tilde{Q}}/\Gamma$,
$S={\tilde{S}}/\Gamma$. The covering projection will be denoted by $p:\tQ\ra\tQ/\Gamma=Q$.
Let $T$ be the normal bundle of $S$, that is, for $z\in S$,
$T_z=\{ v\in T_zQ\, :\, g(v,u)=0,$ for all $u\in T_zS\}\sbs T_zQ$. Write $\pi(v)=z$ if $v\in T_z$, that is,
$\pi:T\ra S$ is the bundle projection. The unit sphere bundle and unit disc bundle of $T$ will be denoted by $N$
and $W$, respectively. Note that the
normal bundle, normal sphere bundle and the normal disc bundle
of $\tS$ in $\tQ$ are the liftings $\tT$, $\tN$ and $\tW$ of $T$, $N$ and $W$,
respectively. 
For $v\in T_qQ$ or $v\in T_q\tQ$, $v\neq 0$, the map $t\mapsto exp_q(tv)$, $t\geq 0$, will be denoted by $c_v$ and its image
will be denoted by the same symbol. Since $\tQ$ is simply connected, $c_v$ is a geodesic ray,
for every $v\in \tN$. We have the following well known facts.
 
 \begin{enumerate}
\item[{1.}] For any closed convex set $C\sbs\tQ$, and a geodesic $c$, the function
$t\mapsto \tid\, (\, c (t)\,,\, C\, )$ is convex.
This implies 2 below.
 
\item[{2.}] Let $c$ be a geodesic ray beginning at some $z\in \tS$. Then either $c\sbs \tS$ or 
${\tilde{d}}\, (\, c(t)\, ,\, \tS\, )\ra \infty$, as $t\ra\infty$.
 
\item[{3.}] For every $v\in T$, $v\neq 0$, $c_v$ is a geodesic ray. Moreover,
for non-zero vectors $v_1,\, v_2\in T$, with $\pi(v_1)\neq\pi(v_2)$, we have that the function
$t\mapsto d\, (\, c_{v_1}(t)\, ,\,c_{v_2}\, )$ tends to $\infty$
as $t\ra\infty$.

\item[{4.}] The exponential map $E:T\ra Q$, $E(v)=exp_{\pi(v)}(v)$, is a diffeomorphism. 
We can define then the submersion $proj:Q\ra S$, $proj(q)=z$, if $exp(v)=q$, for some $v\in T_z$. 
Write also $\eta(q)=|v\, |$ and we have $\eta (q)=d\, (\, q\, ,\, S\, )$.
Also, the exponential map $\tilde{E}:\tT\ra \tQ$, $\tilde{E}(v)=exp_{\pi(v)}(v)$, is a diffeomorphism
and $\tilde{E}$ is a lifting of $E$.

\item[{5.}] Since $S$ is compact there is a map $\varrho\,$ such that: (1) for $q_1,q_2\in Q$, 
$\varrho (a)\, d\, (\, proj(q_1)\, ,\, proj(q_2)\, )\leq
d\, (\, q_1\, ,\, q_2, )$, where $a=min\{ \eta(q_1),\eta(q_2)\}\,\, $ (2) $\varrho (0)=1$,
$\varrho$ is an increasing function and tends to $\infty$ as $t\ra\infty$.

\item[{6.}] Recall that we are assuming that all sectional curvatures of $\tQ$ are less that $c_2<0$.
Given $\lambda\geq 1$, $\epsilon \geq 0$, there is a number $K=K(\lambda,\epsilon,c_2)$ such that the following happens.
For every $(\lambda,\epsilon)$-quasi-geodesic $c$ in $\tQ$ there is a unit speed geodesic $\beta$ with the same endpoints
as $c$, whose Hausdorff distance from $c$ is less or equal $K$. Note $K$ depends on $\lambda ,\epsilon ,c_2$ but not
on the particular manifold $\tQ$ (see, for instance, \cite{BH}, p. 401; see also Proposition 1.2 on p. 399 of \cite{BH}).
\end{enumerate}

Recall that the space at infinity $\bo\tQ$ of $\tQ$ can be defined as $\{\, quasi-geodesic\,\, rays\,\, in
\,\, \tQ \}/\sim$
where the relation $\sim$ is given by: $\beta_1\sim\beta_2$ if their Hausdorff distance is finite.
We say that a quasi-geodesic $\beta$ {\it converges to} $p\in\bo\tQ$ if $\beta\in p$.
Fact 6 implies that we can define $\bo\tQ$ also by
$\{\, geodesics\,\, rays\,\, in\,\,\tQ\, \}/\sim$.
We consider $\bo\tQ$ with the usual cone topology (see \cite{BH}, p. 263). Recall that, for any $q\in\tQ$, the map
$\{ v\in T_q\tQ\, :\, |v|=1\} \ra\bo\tQ$ given by $v\mapsto [c_v]$ is a homeomorphism.
Let $\varsigma :[0,1]\ra[0,\infty)$ be a homeomorphism that is the identity near 0. 
We also have that $\ctQ=\tQ\cup\bo\tQ$ can be given a topology such that 
the map $\{ v\in T_q\tQ\, :\, |v|\leq1\} \ra\bo\tQ$ given by $v\mapsto exp_q (\varsigma  (|v|)\frac{v}{|v|})$,
for $|v|<1$ and $v\mapsto [c_v]$ for ${v}=1$, is a homeomorphism.
We have some more facts or comments.

\begin{enumerate}
\item[{7.}] Given $q\in\tQ$ and $p\in\bo\tQ$ there is a unique unit speed geodesic ray $\beta$ beginning at
$q$ and converging to $p$.

\item[{8.}] Since $\tS$ is convex in $\tQ$ every geodesic ray in $\tS$ is a geodesic ray in $\tQ$. Therefore
$\bo\tS\sbs\bo\tQ$. For a quasi-geodesic ray $\beta$ we have: $[\beta]\in\bo\tQ\sm\bo\tS$ if and only
if $\beta$ diverges from $\tS$, that is ${\tilde{d}}\, (\, \beta(t)\, ,\, \tS\, )\ra \infty$, as $t\ra\infty$.

\item[{9.}] For every $p\in\bo\tQ\sm\bo\tS$ there is a unique $v\in \tN$ such that $c_v$ converges to $p$.
Moreover, the map $\tA:\tN\ra\bo\tQ\sm\bo\tS$, given by $\tA(v)=[c_v]$ is a homeomorphism.
Furthermore, we can extend $\tA$ to a homeomorphism $\tW \ra\ctQ\sm\bo\tS$ by  defining $\tA(v)=\tilde{E}(\varsigma
(|v|)\frac{v}{|v|})=exp_q (\varsigma  (|v|)\frac{v}{|v|})$, for $|v|<1$, $v\in\tW_q$
(recall that $\varsigma$ is the identity near zero). 
\end{enumerate}

\noindent {\bf Lemma 2.2.} {\it Let $\beta:[a,\infty )\ra \tQ$. The following are equivalent.}
\begin{enumerate}
\item[{(i)}] {\it $\beta$ is a quasi-geodesic ray and diverges from $\tilde S$.}

\item[{(ii)}] {\it $p\beta$ is a quasi-geodesic ray.}
\end{enumerate}

{\bf Proof.} First note that if a path $\alpha(t)$, $t\geq a$, satisfies the $(\lambda,\epsilon)$-quasi-geodesic
ray condition, for $t\geq a'\geq a$, then $\alpha(t)$ satisfies the
$(\lambda,\epsilon')$-quasi-geodesic ray condition, for all $t\geq a$, where $\epsilon'=\epsilon+diameter(\alpha([a,a']))$.

(i) implies (ii). Let $\beta$ satisfy (i).
Then there are $\lambda\geq 1$, $\epsilon \geq 0$ such that 
$\frac{1}{\lambda}|t-t'|-\epsilon\leq \tid (\beta (t),\beta(t'))\leq \lambda |t-t'|+\epsilon$,
for every $t,t'\geq a$. 
Fix $t,t'\geq a$ and let $\alpha$ be
the unit speed geodesic segment joining $\beta(t)$ to  $\beta(t')$. Then $p\alpha$ joins 
$p\beta(t)$ to  $p\beta(t')$. Therefore $d (p\beta (t),p\beta (t'))\leq length_{g}(p\alpha)=
length_{\tilde{g}}(\alpha)= d(\beta (t),\beta (t'))\leq \lambda|t-t'|+\epsilon$.
We proved that $d (p\beta (t),p\beta (t'))\leq  \lambda|t-t'|+\epsilon$.

We show the other inequality. By item 6, $\beta$ is at finite Hausdorff distance
(say, $K\geq 0$) from a geodesic ray $\alpha$. Since $\beta$ (hence $\alpha$) gets
far away from $\tS$, it converges to a point at infinity in $\bo\tQ\sm\bo\tS$. Therefore
we can assume that $\alpha(t)=c_{\tv}(t)=exp_{\tilde{z}}(t\tv)$ for some $\tv\in\tT_{\tilde{z}}$, with 
$|\tilde{v}|=1$. It follows that
$p\beta$ is at Hausdorff distance $K'=K+d(\beta (a),\tS)$ from $c_v$, where $v\in T_z$ is the image of
$\tv$ by the derivative $Dp(\tilde{z})$, and $z=p(\tilde{z})$. 
Note that $c_v$ is a geodesic ray in $Q$ (see item 3).
Let $U$ denote the
$K$ neighborhood of $c_v$ in $Q$ and $\tilde{U}$ the $K$ neighborhood of $c_{\tv}$ in $\tQ$.
We claim that $p:\tilde{U}\ra U$ satisfies: $d(p(x),p(y))\geq\tid(x,y)-4K$, for $x,y\in\tilde{U}$.
To prove this let $t,t'\geq$ such that $d(x,c(t))=d(x,c_v)\leq K$ and $d(y,c(t'))=d(y,c_v)\leq K$.
We have 
$\tid(x,y)\leq \tid(x,c_{\tv}(t))+\tid(c_{\tv}(t),c_{\tv}(t'))+\tid(c_{\tv}(t'),y)\leq 2K+|t-t'|=2K+d(c_v(t),c_v(t'))\leq
2K+ d(c_v(t),p(x))+d(p(x),p(y))+d(p(y),c_v(t'))\leq 4K+d(p(x),p(y))$. This proves our claim.
Consequently $d(p\beta(t),p\beta(t'))\geq\tid(\beta(t),\beta(t'))-4K\geq\frac{1}{\lambda}|t-t'|-(\epsilon+4K)$.

(ii) implies (i). Let $\beta$ satisfy (ii). Since $p\beta$ is a proper map its
distance to $S$ must tend to infinity. Hence the distance of $\beta$ to $\tS$
also tends to infinity. 

Let $p\beta$ satisfy
$\frac{1}{\lambda}|t-t'|-\epsilon\leq d (p\beta (t),p\beta(t'))\leq \lambda |t-t'|+\epsilon$,
for some $\lambda\geq 1$, $\epsilon \geq 0$. Fix $t,t'\geq a$ and let $\alpha$ be
the unit speed geodesic segment joining $\beta(t)$ to  $\beta(t')$. Then $p\alpha$ joins 
$p\beta(t)$ to  $p\beta(t')$. Therefore $\tid (\beta (t),\beta (t'))=length_{\tilde{g}}(\alpha)=
length_g(p\alpha)\geq d(p\beta (t),p\beta (t'))\geq \frac{1}{\lambda}|t-t'|-\epsilon$.
It follows that $\frac{1}{\lambda}|t-t'|-\epsilon\leq \tid (\beta (t),\beta(t'))$.

We prove the other inequality.
Since $S$ is compact and by item 5, the radius of injectivity of a point in $Q$ tends
to infinity as the points gets far from $S$. Hence there is $a'\geq a$ such that
for every $t\geq a'$, the ball of radius $e=\lambda+\epsilon$ centered
at $\beta(t)$ is convex. Let $t'>t>a'$ and $n$ an integer such that
$n< t'-t\leq n+1$. Let $\alpha_k$, $k=1,...,n$, be the unit speed geodesic segment from $p\beta(t+k-1)$ to 
$p\beta (t+k)$, and $\alpha_{n+1}$ the unit speed geodesic segment from $p\beta(t+n)$ to 
$p\beta (t')$. Note that $length_{g}(\alpha_k)=d(p\beta(t+k-1),p\beta (t+k))\leq \lambda+\epsilon=e$. 
Therefore $p\beta|_{[t+k-1,t+k]}$ is homotopic, rel endpoints, to $\alpha_k$ (analogously for $\alpha_{n+1}$).
Let $\alpha$ the concatenation $\alpha_1*...*\alpha_{n+1}$. Then $\alpha$ is homotopic, rel endpoints,
to $p\beta|_{[t,t']}$. Note that the length of $\alpha$ is $\leq (n+1)e$.
Let $\tilde{\alpha}$ be the lifting of $\alpha$ beginning at $\beta(a')$. Then
$\tilde{\alpha}$ is homotopic, rel endpoints, to $\beta|_{[t,t']}$. Hence $\tid (\beta(t),\beta(t'))
\leq length (\tilde{\alpha})\leq (n+1)e=ne+e<e(t'-t)+e$. We showed that
$\frac{1}{\lambda}|t-t'|-\epsilon\leq\tid (\beta(t),\beta(t'))<(\lambda+\epsilon)|t'-t|+(\lambda+\epsilon)$. 
This proves the Lemma.
\vspace{.2in}

Let $Q_1$, $Q_2$ be two complete simply connected negatively curved manifolds.
If $\beta$ is a quasi-geodesic in $Q_1$ and $f:Q_1\ra Q_2$  is a quasi-isometry then
$f(\beta)$ is also a quasi-geodesic. Also, if two subsets of $Q_1$ have finite Hausdorff
distance, their images under $f$ will have finite Hausdorff distance as well. Therefore $f$ induces
a map $f_\infty:\bo Q_1\ra\bo Q_2$. Hence $f$ extends to $\bar{f}:\cQ_1\ra\cQ_2$ by $\bar{f}|_{\bo Q_1}=f_\infty$
and $\bar{f}|_{Q_1} =f$.
We have

\begin{enumerate}
\item[10.] For every quasi-isometry $f:Q_1\ra Q_2$,  $f_\infty:\bo Q_1\ra\bo Q_2$ is a homeomorphism. 
In addition, if $f$ is a homeomorphism, then $\bar{f}$ is a homeomorphism.

\item[{11.}] Let $g'$ be another complete Riemannian metric on $\tQ$ 
whose sectional curvatures are also $\leq c_2 <0$, and
such that there are constants $a,b>0$ with $a^2\leq g'(v,v)\leq b^2$ for every $v\in T\tQ$ with $\tilde{g}(v,v)=1$,
and such that $\tS$ is also a convex subset of $(\tQ,g')$.
Then $\bo\tQ$ is the same if defined using $\tilde{g}$ or $g'$.
Moreover item 9 above also holds for $(\tQ,g')$ (with respect to all proper concepts defined using $g'$ instead
of $\tilde{g}$). This is because the identity $(\tQ,\tilde{g})\ra (\tQ,g')$ induces the  homeomorphism 
$\bo\tQ\ra\bo\tQ$ that preserves $\bo\tS$ (see Lemma 2.1 and item 10).
\end{enumerate}

Since $\Gamma$ acts by isometries on $\tQ$, we have that $\Gamma$ acts on $\bo\tQ$ (see item 10).
Also, since $\Gamma$ preserves $\tS$, $\Gamma$ also preserves $\bo\tS$. Hence $\Gamma$ acts
on $\bo\tQ\sm\bo\tS$. Since $S$ is closed, we have

\begin{enumerate}
\item[{12.}] For every $\gamma\in\Gamma$,  $\gamma :\bo\tQ\sm\bo\tS\ra\bo\tQ\sm\bo\tS$ has no fixed points.
Therefore the action of $\Gamma$ on
$\ctQ\sm\bo\tS$ is free. Moreover, the action of $\Gamma$ on $\ctQ\sm\bo\tS$ is properly
discontinuous.
\end{enumerate}

We now define the space at infinity $\bo Q$ of $Q$ as $\{\, quasi-geodesic\,\, rays\,\, in\,\, Q \}/\sim$.
As before, the relation $\sim$ is given by: $\beta_1\sim\beta_2$ if their Hausdorff distance is finite.
We can define a topology on $\bo Q$ in the same way as for $\bo\tQ$, but we can take advantage of
the already defined topology of  $\bo\tQ$.\\

\noindent {\bf Lemma 2.3.} {\it There is a one-to-one correspondence between $\bo Q$ and 
$\left(\bo\tQ\sm\bo\tS\right)\, /\,\Gamma$.}\vspace{.1in}

\noindent {\bf Proof.} By path lifting and Lemma 2.2 there is a one-to-one correspondence between the
sets $\{ \, quasi-geodesic\,\, rays\,\, in\,\, Q \}\,$ and 
$\,\{\, quasi-geodesic\,\, rays\,\, in\,\, \tQ\,\, that \,\, diverge\,\, from\,\,\tS \}\, /\, \Gamma$.
Then the correspondence $[\beta]\mapsto p\, (\beta)$, for quasi-geodesic rays in $\tQ$ that diverge from $\tS$
is one-to-one (see item 8).
This proves the Lemma.\\

We define then the topology of $\bo Q$ such that the one-to-one correspondence
mentioned in the proof of the Lemma is a homeomorphism.
Also, we define the topology on $\cQ= Q\cup\bo Q$ such that $\left(\ctQ\sm\bo\tS\right)\, /\, \Gamma\ra\cQ$
is a homeomorphism. It is straightforward to verify that $Q$ and $\bo Q$ are subspaces of $\cQ$
(see also item 12). The next Lemma is a version of item 9 for $Q$.\\

\noindent {\bf Lemma 2.4.} {\it
For every $p\in\bo Q$ there is a unique $v\in N$ such that $c_v$ converges to $p$.
Moreover, the map $A:N\ra\bo Q$, given by $A(v)=[c_v]$ is a homeomorphism.
Furthermore, we can extend $A$ to a homeomorphism $W \ra\bo Q$ by  defining 
$A(v) = E( (\varsigma  (|v|)\frac{v}{|v|}))$, for $|v|<1$. (Recall $\varsigma$ is the identity near 0.)
Also, $\tA$ is a lifting of $A$.}\vspace{.1in}

\noindent {\bf Proof.} The first statement follows from items 4 and 5.
Define $A(v)=p\tA(\tilde{v})$, where $Dp(\tilde{v})=v$. Items 9 and 12 imply
the Lemma. See also item 4.\\

\noindent We will write $\eta ([c_v])=\infty$ and  $E(\infty v)=[c_v]$, for $v\in N$
(see item 5).\\

\noindent {\bf Lemma 2.5.} {\it Let $v\in N$ and $q_n=E(t_nv_n)$, $t_n\in [0,\infty]$, $v_n\in T$ 
and $|v_n|$ bounded away  from both $0$ and $+\infty$. Then 
$q_n\ra [c_v]$ (in $\bo Q$) if and only if $t_n\ra\infty$ and $v_n\ra v$ .}\vspace{.1in}

\noindent {\bf Proof.} It follows from Lemma 2.4.\\

\noindent We also have a version of item 11 for $Q$.\\

\noindent {\bf Lemma 2.6.} {\it Let $g'$ be another complete Riemannian metric on $Q$ 
whose sectional curvatures are also $\leq c_2 <0$, and
such that 
there are constants $a,b>0$ with $a^2\leq g'(v,v)\leq b^2$ for every $v\in TQ$ with $g(v,v)=1$,
and such that $S$ is also a convex subset of $(Q,g')$.
Then $\bo Q$ is the same if defined using $g$ or $g'$.
Moreover Lemma 2.4 and 2.5 above also holds for $(Q,g')$ (with respect to all proper concepts defined using $g'$ instead
of $g$).}\vspace{.1in}

\noindent {\bf Proof.} It follows from item 11 and Lemma 2.5.
Note that the liftings $\tilde{g}$, $\tilde{g}'$ of $g$ and $g'$, satisfy 
$a^2\leq \tilde{g}'(v,v)\leq b^2$ for every $v\in T\tQ$ with $\tilde{g}(v,v)=1$.
This proves the Lemma.
\vspace{.6in}

\noindent {\bf \large  Section 3. Proof of Theorem 1.}\\

Let the metric $g$ and the closed simple curve $\alpha$ be as in the statement of the Theorem.
Write $N=\bS^1\times\bS^{n-2}$ and $\Sigma^M=\Lambda_g\Phi^M$, where $\Phi^M=\Phi^M(\alpha,V,r)$.
The base point of the $k$-sphere $\bS^k$ will always
be the point $u_0=(1,0,...,0)$.
Let $\theta :\bS^k\ra DIFF(N\times I, \p)$, $\theta (u_0)=1_{N\times I}$, represent an element in $\pi_k (\,
DIFF(N\times I,\p )\, )$. \\

We will prove that if $\pi_k(\Sigma^M)([\theta])$ is zero, then $\pi_k (\iota_N)([\theta])$ is also zero.
Equivalently, if $\Sigma^M\, \theta$ extends to the ($k+1$)-disc $\D^{k+1}$, then $\iota_N\theta$ also extends
to $\D^{k+1}$.
So, suppose that $\Sigma^M\,\theta :\bS^k\ra\mo (M)\,$ extends to a map $\sigma' :\D^{k+1}\ra\mo (M)$.
We can assume that this map is smooth.\\

\noindent {\bf Remark.}  Originally $\sigma'$ may not be smooth, but it is homotopic to a smooth map.
By ``$\sigma'$ is smooth" we mean that the map $\D^{k+1}\times (TM\oplus TM)\ra \R$,
given by $(u, v_1,v_2)\mapsto \sigma'(u)_x(v_1,v_2)$, $v_1,v_2\in T_xM$, is smooth.
To homotope a given $\sigma'$ to a smooth one $\sigma''$ we can use classical averaging techniques:
just define $\sigma_x(u)''(v_1,v_2)=\int_{\R^{k+1}}\eta(u-w)\,\sigma'(w)_x(v_1,v_2)\, dw$, which is smooth. Here:
 (1)\,  $\eta$ is a smooth $\epsilon$-bump function, i.e. $\int_{\R^{k+1}}\eta =1$ and $\eta(w)=0$, for $|w|\geq\epsilon$
and, (2) we are extending $\sigma'$ (originally defined on $\D^{k+1}$) to all $\R^n$, radially.
Since $\sigma'$ is continuous, the second order derivatives of $\sigma'_x(u)$ and $\sigma'_x(u')$ are close
for $u$ close to $u'$. Therefore the second order derivatives of $\sigma'_x(u)$ are close
to the second order derivatives of $\sigma''_x(u)$. Hence, if $\epsilon$ is sufficiently small,
we will also have $\sigma''(u)\in\mo (M)$.\\

Also, by deforming $\sigma'$, we can assume that it is radial near $\p\,\D^{k+1}$.
Thus $\sigma '(u),\,\, u\in\D^{k+1},$ is a negatively curved metric on $M$. Also, $\sigma '(u)=\Sigma^M\,\theta (u)$,
for $u\in\bS^k$, and $\sigma'(u_0)=g$. Since $\sigma'$ is continuous there is a constant $c_2<0$ such that
all sectional curvatures of the Riemannian manifolds $(M,\sigma'(u))$, $u\in\D^{k+1}$, are less or equal $c_2$.
Write $\varphi_u=\Phi^M(\theta(u))$, $u\in\bS^k$.
Hence we have that $\sigma'(u)=(\varphi_u)_*\sigma'(u_0)=(\varphi_u)_*g$, for $u\in\bS^k$.
Note that $\varphi_u$ is, by definition, 
the identity outside the closed normal geodesic tubular neighborhood $U$ of width $2r$ of $\alpha$.
Also, $\varphi_u$ is the identity on the closed normal geodesic tubular neighborhood of width $r$ of $\alpha$.
Note that $\varphi_u:M\ra M$ induces the identity at the $\pi_1$-level
and hence $\varphi_u$  is freely homotopic to $1_M$.\\

Since $\sigma'$ is continuous and $\D^{k+1}$ is compact we can find constants
$a,b>0$ such that $a^2\leq \sigma'(u)(v,v)\leq b^2$ for every $v\in TM$ with $g(v,v)=1$,
$u\in\D^{k+1}$.\\

Let $Q$ be the covering space of $M$ with respect to the infinite cyclic subgroup
of $\pi_1 (M,\alpha (1))$ generated by $\alpha$. Denote by $\sigma(u)$ the pullback on $Q$
of the metric $\sigma'(u)$ on $M$. For the lifting of $g$ on $Q$ we use the same letter $g$.
Note that $\alpha$ lifts to $Q$ and we denote this lifting
also by $\alpha$.
Let $\phi_u:Q\ra Q$ be diffeomorphism which is the unique lifting of $\varphi_u$ to $Q$
with the property that $\phi_u|_{\alpha}$ is the identity.
We have some comments.
\begin{enumerate}
\item[{(i.)}] $\sigma(u)=(\phi_u)_*\sigma(u_0)=(\phi_u)_*g$, for $u\in\bS^k$.

\item[{(ii.)}] The tubular neighborhood $U$ lifts to a countable number of components, with
exactly one being diffeomorphic to $U$. We call this lifting also by $U$. All
other components $U_1, \, U_2,...$ are diffeomorphic to $\D^{n-1}\times\R$.
Note that $\phi_u$ is the identity outside the union of $\bigcup U_i$ and $U$
and inside the closed normal geodesic tubular neighborhood of width $r$ of $\alpha$.

\item[{(iii.)}] Since $\varphi_u:M\ra M$ induces the identity at the $\pi_1$-level,
and $\bS^{k}$ is compact, there is a constant $C$ such that $d_{\sigma(u')}(\, p\, ,\,\phi_u(p)\, )<C$,
for any $u,u'\in\bS^k$, where $d_{\sigma(u')}$ denotes the distance in the Riemannian
manifold $(Q,\sigma(u'))$.

\item[{(iv.)}] $(\phi_u)|_U=\left[\Phi^Q(\alpha, V',r)\theta (u)\right]|_U$, for $u\in\bS^k$. Here $V'$ is the lifting of $V$.

\item[{(v).}] We have that $a^2\leq \sigma(u)(v,v)\leq b^2$ for every $v\in TQ$ with $g(v,v)=1$,
$u\in\D^{k+1}$. It follows that $\frac{a^2}{b^2}\leq \sigma(u)(v,v)\leq \frac{b^2}{a^2}$ for every $v\in TQ$ with 
$\sigma(u')(v,v)=1$, $u,u'\in\D^{k+1}$.

\item[{(vi.)}] All sectional curvatures of the Riemannian manifolds $(Q,\sigma(u))$, $u\in\D^{k+1}$, are less or equal $c_2$.
\end{enumerate}

Since $(M,\sigma'(u))$ is a closed negatively curved manifold, it
contains exactly one immersed closed geodesic freely homotopic to $\alpha\sbs M$.
Therefore $(Q,\sigma(u))$ contains exactly one embedded closed geodesic $\alpha_u$ freely homotopic to $\alpha\sbs Q$.
Note that $\alpha_u$ is unique up to affine reparametrizations. 
Write $\alpha_0=\alpha_{u_0}$ and note that $\alpha_u =\phi_u(\alpha_0)$, for all $u\in\bS^k$.\\

Since $n\geq 5$ we can find a compactly supported smooth isotopy $s:Q\times I\ra Q$, $s_0=1_Q$, with
$s_1(\alpha_0)=\alpha$. 
Using $s$ we get a homotopy $(s_t)^{-1} \phi_u s_t$ between $\phi_u$ and $\psi_u=(s_1)^{-1}\phi_u s_1$.
Therefore we can assume that for $u\in\bS^k$ we have $\sigma(u)=(\psi_u)_*g$.
Note that (ii) above still holds with $U'=(s_1)^{-1}U$, $U'_i=(s_1)^{-1}U_i$ instead of $U$, $U_i$, respectively.
Note that $U'_i$ coincides with $U_i$ outside a compact set.
Also, since $s$ is compactly supported (iii) holds too. For (iv) we assume that $U'$ is the
closed normal geodesic tubular neighborhood of width $2r$ of $\alpha_0$ and $s_1$ sends geodesic of length
$2r$ beginning orthogonally at $\alpha_0$ isometrically to
geodesic of length $2r$ beginning orthogonally at $\alpha$ (we may have to consider a much smaller $r>0$ here). 
Note that (v) and (vi) still hold.
The following version of (iv) is true
\begin{enumerate}
\item[{(iv'.)}] $(\psi_u)|_{U'}=\left[\Phi^Q(\alpha_0, V'',r)\theta (u)\right]|_{U'}$, for $u\in\bS^k$. 
Here $V''=(s_1^{-1})_*V'$.
\end{enumerate}

Now, by [6, Prop. 5.5] $\alpha_u$ depends smoothly on $u\in\D^{k+1}$. Hence
we have a smooth map $h:\D^{k+1}\times\bS^1\ra Q $, given by $h_u=\alpha_u$.
Note that $h$ is radial near $\p$.
We have the following facts:
\begin{enumerate}
\item[{1.}] We can identify $\bS^1$ with its image $\alpha_0$ and, using the exponential
map orthogonal to $\bS^1$, with respect to $g=\sigma(u_0)$ and the trivialization $V''$, we can identify $Q$ to 
$\bS^1\times\R^{n-1}$. With this identification $V''$ becomes just the canonical base $E=\{ e_1,...,e_{n-1}\}$
and (iv') above has now the following form:
$(\psi_u)|_{U'}=\left[\Phi^Q(\alpha_0, E,r)\theta (u)\right]|_{U'}$, for $u\in\bS^k$. 

\item[{2.}]  Because of the argument above (using the homotopy $s$) we can not guarantee that
all metrics $\sigma(u)$ are lifted metrics from $M$, but we do have that all liftings of the $\sigma (u)$
to the universal cover ${\tilde Q}={\tilde M}$ are all quasi-isometric.
\end{enumerate}

The next Claim says that we can assume all $h_u=\alpha_u:\bS^1\ra Q$
to be equal to $\alpha_0$.\\

\noindent {\bf Claim 1.} {\it We can modify $\sigma$ (hence also $\alpha_u$ and $h$)
on $int\, (\D^{k+1})$ such that }
\begin{enumerate}
\item[{a.}] {\it The liftings of the metrics $\sigma(u)$ to the universal cover
${\tilde Q}={\tilde M}$ are all quasi-isometric.}

\item[{b.}] {\it $\alpha_u=\alpha_0$, for all $u\in\D^{k+1}$.}
\end{enumerate}

\noindent {\bf Proof of Claim 1.} 
Let $H$ be as in Lemma 1.4. 
Then the required new metrics are just $[(H_u)_1]^*\sigma(u)$,
that is, the pull-backs  of $\sigma(u)$ by the inverse of the diffeomorphism
given by the isotopy $H_u$ at time $t=0$. 
Note that the metrics do not change outside a compact set of $Q$.
Just one more detail. In order to be able to apply Lemma 1.4 for $k=0$ we have to know
that the loop $\beta: \D^1\ra Q$ given by $\beta (u)= h(u,1)$ is homotopy trivial. But if this is not the case
let $l$ be such that $\beta$ is homotopic (rel base point) to $\alpha_0^{-l}$. Then just replace
$h$ by $h\, \vartheta$, where $\vartheta:\D^1\times\bS^1\ra \D^1\times\bS^1$, $\vartheta (u,z)=(u,e^{\pi l(u+1)i}.\, z)$.
Note that $h_u$ and $(h\,\vartheta)_u$ represent the same geodesic, but with different basepoint.
This proves Claim 1.\\

Hence, from now on, we assume that all $\alpha_u$ are equal to $\alpha_0:\bS^1\ra Q$.
Note that the new metrics $\sigma(u)$, $u\in int\, (\D^{k+1})$, are not necessarily
pull-back from metrics in $M$.
Recall that we are identifying $Q$ with $\bS^1\times\R^{n-1}$, and the rays $\{ z\}\times \R^+ v$, 
$v\in \bS^{n-2}$, are geodesics (with respect to $g=\sigma(u_0) $) emanating from $z\in\bS^1\sbs Q$  
and normal to $\bS^1$.
Denote by $W_\delta=\bS^1\times\D^{n-1}(\delta )$ the closed normal tubular neighborhood of $\bS^1$ in $Q$
of width $\delta>0$, with
respect to the metric $\sigma(u_0)$. Note that $\p W_\delta=\bS^1\times\bS^{n-2}(\delta)$.\\

For each $u\in\D^{k+1}$ and $z\in\bS^1$, let $T^u(z)$ be the orthogonal complement of the
tangent space $T_z\bS^1\sbs T_zQ$ with respect to the $\sigma(u)$ metric
and denote by $exp^u_z:T^u(z)\ra Q $ the normal exponential map, also with respect to the
$\sigma(u)$ metric. Note that the map $exp^u:T^u\ra Q$ is a diffeomorphism,
where $T^u$ is the bundle over $\bS^1$ whose fibers are $T^u(z)$, $z\in\bS^1$.
We will denote by $N^u$ the sphere bundle of $T^u$.
The orthogonal projection (with respect to the $\sigma(u_0)$ metric)
of the tangent vectors $(z, e_1),...,(z,e_{n-1})\in T_z Q=\{ z\}\times\R^{n-1}$
(here $e_1=(1,0,...,0)$, $e_2=(0,1,0,...,0)$,... ) into $T^u(z)$ gives a base of
$T^u(z)$. Applying the Gram-Schimidt orthogonalization process we obtain and orthonormal
base $v_u^1(z),...,v_u^{n-1}(z)$ of $T^u(z)$. Clearly, these bases are continuous in $z$, hence
they provide a trivialization of the normal bundle $T^u$. 
We denote by $\chi_u:T^u\ra \bS^1\times\R^{n-1}$ the bundle 
trivializations given by $\chi_u(v_u^i(z))=(z,e_i)$.
Note that these trivializations are continuous in $u\in \D^{k+1}$.\\

For every $(u,z,v)\in\D^{k+1}\times\bS^1\times (\R^{n-1}\sm\{ 0\})$ define 
$\tau_u(z,v)=(z',v')$, where $\chi_u \circ (exp^u)^{-1}(z,v)=(z', w)$ and $v'=\frac{w}{|w|}$.
Then $\tau_u: \bS^1\times (\R^{n-1}\sm\{ 0\}) \ra \bS^1\times\bS^{n-2}$ is a smooth map. The restriction of
$\tau_u$ to any $\p W_\delta\sbs\bS^1\times\R^{n-1}$ will be denoted also by $\tau_u$.
>From now on we assume $\delta<r$.\\

\noindent {\bf Claim 2.} {\it There is $\delta >0$ such that
the map $\tau_u: \p W_\delta \ra \bS^1\times\bS^{n-2}$ is a diffeomorphism.}\vspace{.1in}

\noindent {\bf Proof of Claim 2.} 
Just apply Lemma 1.6 to the map $\chi_u \circ (exp^u)^{-1} $.
This proves Claim 2.\\

Note that $\tau_u $ depends continuously on $u$. Note also that Claim 2 implies that every normal geodesic 
(with respect to any metric $\sigma(u)$) emanating from 
$\alpha_0$,  intersects $\p W_\delta$ transversally in a unique point.
Denote by $\rho_u :\p W_\delta\ra (0,\infty )$
the smooth map given by $\tau_u(z,v)=|w|$, where we are using the notation before the statement of Claim 2.\\

To simplify our notation we take $\delta =1$ and write $W=W_1$. Thus $\p W=N=\bS^1\times\bS^{n-2}$
and we write $ N\times [1,\infty )=Q\sm int\, W$.
Now, for each $u\in\D^{k+1}$ we define a self-diffeomorphism $f_u\in DIFF(N\times [1,\infty ), N\times\{ 1\})$
by $$f_u((z,v),t)=exp^u_{z'}(\,\,  [\chi_u]^{-1}(\, \, z',\, \rho_u(z,v)\, tv'\, )\,\, )$$

\noindent where $\tau_u(z,v)=(z',v')$. It is not difficult to show that $f_u((z,v),1)=((z,v),1)$ 
and that $f_u$ is continuous in $u\in\D^{K+1}$.\\

Here is an alternative interpretation of $f_u$. For $(u,z,v)\in \D^{k+1}\times\bS^1\times T^u(z)$, denote by
$c^u_{(z,v)}:[0,\infty )\ra Q$ the $\sigma(u)$ geodesic ray given by
$c^u_{(z,v)}(t)=exp^u_z(tv)$. Then $f_u$ sends $c^{u_0}_{(z,v)}$ to  $c^u_{(z',s)}$, where
$exp^u_{z'}(s)=(z,v)\in Q$. Explicitly, we have $f_u (\,  c^{u_0}_{(z,v)}(t)   \, )= c^u_{(z',s)}(|s|t)$, for 
$t\geq 1$.  Using Claim 2 it is not difficult to prove that $f_u ( \, N\times [1,\infty ) \, )=N\times [1,\infty )$ 
and that $f_u$ is a diffeomorphism.\\

We denote by $\p_\infty Q$ the space at infinity of $Q$ with respect to the $\sigma(u_0)$ metric.
Recall that the elements of $\p_\infty Q$ are equivalence classes $[\beta]$ of $\sigma(u_0)$ 
quasi-geodesic rays $\beta :[a,\infty)\ra Q=\bS^1\times\R^{n-1}$ (see Section 2). Note that, since all 
metrics $\sigma(u)$ are quasi-isometric, a 
$\sigma(u)$ quasi-geodesic ray  is a $\sigma(u')$ quasi-geodesic ray, 
for any $u,u'\in\D^{k+1}$. Hence
$\p_\infty Q$ is independent of the metric $\sigma(u)$ used  (see (v) and Lemma 2.?). Still, the choice of a $u\in\D^{k+1}$,
gives canonical elements in each equivalence class in $\p_\infty Q$: just choose the unique unit speed 
$\sigma(u)$ geodesic ray
that ``converges'' (that is,``belongs'') to the class, and that emanates
$\sigma(u)$-orthogonally from $\bS^1\sbs Q$. If we choose the $\sigma(u_0)$ metric, this set of geodesic
rays is in one-to-one correspondence with $N=\bS^1\times\bS^{n-2}\sbs Q$. 
We identify $N\times \{\infty\}$ with $\p_\infty Q$ by: $((z,v),\infty )\mapsto [c^{u_0}_{(z,v)}]$.
Hence we can write now $(Q\sm int\, W\, )\,\cup\,\p_\infty Q
=(\, N\times [1,\infty )\, )\,\cup\,\p_\infty Q=N\times [1,\infty ]$ (see Lemma 2.?).\\

We now extend each $f_u$ to a map $f_u:N\times[1,\infty ]\ra N\times[1,\infty ]$
in the following way. For $((z,v),\infty )=[c^{u_0}_{(z,v)}]$ define
$f_u( \, [c^{u_0}_{(z,v)}\, ] \,  )\, =\, [\, f_u( \, c^{u_0}_{(z,v)}\, )\, ]$.
Recall that, as we mentioned before, 
we have $f_u (\,  c^{u_0}_{(z,v)}(t)   \, )= c^u_{(z',s)}(|s|t)$, for 
$exp^u_{z'}(s)=(z,v)\in Q$, $t\geq 1$. That is, $f_u( \, c^{u_0}_{(z,v)}\, )$
is a $\sigma(u)$ geodesic ray, hence it is a $\sigma(u_0)$ quasi-geodesic ray.
Therefore $[\, f_u( \, c^{u_0}_{(z,v)}\, )\, ]$ is a well defined element in 
$\p_\infty$. \\

\noindent We will write $exp=exp^{u_0}$. Also,
as in Section 2, we will write $exp\, (\infty v)=[c_v]$, for $v\in N$.\\

\noindent {\bf Claim 3.} {\it  $f_u:N\times [1,\infty]\ra N\times [1,\infty]$ is a homeomorphism.}\vspace{.1in}

\noindent {\bf Proof of Claim 3.}
Note that $f_u$ is already continuous (even differentiable) on $Q$.
We have to prove that $f_u$ is continuous on points in $\bo Q$.
Let $q_n=exp(t_nv_n)\ra [c_v]$, $v,v_n\in N$, $t_n\in[0,\infty]$. Then, by Lemma 2.5,
$v_n\ra v$ and $t_n\ra\infty$.  Let $u\in\D^{k+1}$ and write $f=f_u$. We have to prove that $q'_n=f(q_n)$ converges to
$f([c_v])=[f(c_v)]$. Write $w_n=(exp^u)^{-1}(v_n)$. Then $w_n\ra w=(exp^u)^{-1}(v)\neq 0$.
Note that $f([c_v])=[f(c_v)]=[c^u_w]$, where $c^u_w$ is the $\sigma (u)$ geodesic ray
$t\mapsto exp^u(tw)$. Note also that, by definition, $f(q_n)=exp^u(t_nw_n)$.
The Claim follows now from Lemmas 2.5 and 2.6.\\

\noindent {\bf Claim 4.} {\it $f_u$ is continuous in $u\in\D^{k+1}$.}\vspace{.1in}

\noindent {\bf Proof of Claim 4.}
Note that we know that $u\mapsto f_u|_{Q}$ is continuous.
Let $q_n=exp(t_nv_n)\ra [c_v]$, $v,v_n\in N$, $t_n\in[0,\infty]$. Then, by Lemma 2.5,
$v_n\ra v$ and $t_n\ra\infty$. Let also $u,u_n\in\D^{k+1}$ with $u_n\ra u$.
To simplify our notation we assume that $u=u_0$ (the proof for a general $u$ is obtained
by properly writing the superscript $u$ on some symbols; see also Lemma 2.6).
Hence, by the previous identifications, $exp^{u_0}=exp:T=Q\ra Q$ is just the identity and
$f_{u_0}$ is also the identity .
Write $f_n=f_{u_n}$ and $w_n=(exp^{u_n})^{-1}(v_n)$. Then 
$w_n\ra (exp^{u_0})^{-1}(v)=v$.
We have to prove that $q'_n=f_n(q_n)=exp^{u_n}(t_nw_n)=c^{u_n}_{w_n}(t_n)$ converges to $f([c_v])=[c_v]$. 
Note that $c^{u_n}_{w_n}(1)=exp^{u_n}(w_n)=v_n\ra v$.
To prove that $q'_n\ra[c_v]$ we will work in $\tQ$ instead of $Q$. Therefore we ``lift''
everything to $\tQ$ and we express this by writing the superscript {\it tilde} over
each symbol. Hence we have $\tv, \tw_n\in\tN$, $u,u_n\in\D^{k+1}$, $t_n >0$ satisfying
\begin{enumerate}
\item[{1.}] $\tw_n\ra\tv$ and $c^{u_n}_{\tw_n}(1)=exp^{u_n}(\tw_n)\ra \tv$.

\item[{2.}] $u_n\ra u_0$, hence $\tilde{\sigma} (u_n)\ra\tilde{\sigma} (u_0)=\tilde{g}.$
\end{enumerate}

We have then that $c_{\tv}$ is a $\tilde{g}$ geodesic ray and the $c^{u_n}_{\tw_n}$ are
$\tilde{\sigma}(u)$ geodesic rays. Write $c^n=c^{u_n}_{\tw_n}$ and $\tilde{q}_n'=c^n(t_n)$.
We have to prove that $\tilde{q}'_n\ra [c_{\tv}]$.
Since $u_n\ra u_0$, the maps $exp^{u_n}\ra exp=1_{\tQ}$ (in the compact-open topology). Therefore
\begin{enumerate}
\item[{(*)}]{\it  for any $R, \, \delta>0$ there is $n_0$ such that $\tid ( c^n (t),c_{\tv} (t))<\delta$, for
$t\leq R$, and $n\geq n_0$.}
\end{enumerate}

Since $c_{\tilde{v}}$ is a unit speed geodesic (i.e. a (1,0)-quasi-geodesic ray), by (1) and (2), for large $n$ we
have that $c^n=c^{u_n}_{\tw_n}$ is a $\tilde{\sigma}(u)$ $(2,0)$-quasi-geodesic ray.
By (v) above and Lemma 2.1 the identity $(\tQ,\tilde{\sigma}(u))\ra (\tQ,\tilde{g})$ is a
$(\lambda,0)$-quasi isometry, where $\lambda=max\{\frac{a^2}{b^2},\frac{b^2}{a^2}\}$.
Therefore, we have that $c^n$ is a $\tilde{g}$ $(2\lambda,0)$-quasi-geodesic ray.
Let $K=K(2\lambda,0,c_2)$ be as in item 6 of Section 2, and $c_2$ is as in (vi) above.
Then there is a unit speed $\tilde{g}$ geodesic ray $\beta_n(t)$, $t\in [1,a_n ]$, that is at $K$ Hausdorff distance from
$c^n$, $t\in [1,t_n]$, and has the same endpoints: $\beta_n(1)=c^n(1)\ra\tv$ and $\beta_n(a_n)=c^n(t_n)=\tilde{q}'_n$.  
Note that $a_n\ra \infty$ because $t_n\ra \infty$.
We have that (*) above \,(take $\delta =1$ in (*))\,
imply that

\begin{enumerate}
\item[{(**)}]{\it given an $R>0$ there is a $n_0$ such that 
$\tid (c_{\tv}(t),\beta_n)\leq C=K+1$, for $t\leq R$ and $n\geq n_0$.}
\end{enumerate}

Since $\tQ$ is complete and simply connected, we can extend each $\beta_n$ to a geodesic
ray $\beta_n:[1,\infty]\ra\tQ$.  Then $[\beta_n]\in\bo\tQ$.
Let $\beta'_n(t)$, $t\in[1,\infty]$ be the unit speed $\tilde{g}$ geodesic ray
with $\beta'_n(1)=\tv$, $\beta'_n(\infty)=\beta_n(\infty)$. Therefore $\tid (\beta_n(t),\beta'_n(t))\leq
\tid (\beta_n(1),\beta'_n(1))=\tid (c^n(1),\tv)\ra 0$. We can assume then that 
$\tid (\beta_n(t),\beta'_n(t))\leq 1$, for all $n$ and $t\geq 1$. Hence, a version of (**) holds 
with $\beta'_n$ instead of $\beta_n$ and $C+1$ instead of $C$.
This new version of (**) implies that $[\beta'_n]\ra [c_{\tilde{v}}]$, and this
together with condition (1.) imply $\beta'_n(t)\ra c_{\tv}(t)$,
for every $t\in [1,\infty]$.  Since $[\beta'_n]\ra [c_{\tv}]$ and $a_n\ra\infty$,
we have that $\beta'_n(a_n)\ra [c_{\tv}]$. But $\tid (\tilde{q}'_n,\beta'_n(a_n))= 
\tid (\beta_n(a_n),\beta'_n(a_n))\leq 1$, therefore $\tilde{q}'_n\ra [c_{\tv}]$. This proves the Claim.
\vspace{.2in}

\noindent {\bf Claim 5.} {\it For all $u\in\bS^k$ we have $f_u|_{Q\sm W}=(\psi_u)|_{Q\sm W}$ 
and $(f_u)|_{\p_\infty}=1_{\p_\infty}$.}\vspace{.1in}

\noindent {\bf Proof of Claim 5.} 
Let $u\in\bS^k$. Since $\sigma(u)=g$ on $W$, then $T^u=T^{u_0}=\bS^1\times\R^{n-1}$
and $exp^u_z(v)=(z,v)$ for all $z\in\bS^1$ and $|v|\leq  1$. It follows that
$f_u (\,  c_{u_0}(z,v)(t)   \, )= c_u(z,v)(t)$, for 
$t\geq 1$.  On the other hand, since $\sigma(u)=(\phi_u)_*\sigma(u_0)$ we have that $\psi:
(Q,\sigma(u_0))\ra (Q,\sigma(u))$ is an isometry. Hence $\psi_u (c_{u_0}(z,v)(t))$, $t\geq 0$, is 
a $\sigma(u)$ geodesic. Since $\psi_u$ is the identity in $W\sbs U'$, we have $\psi_u(z)=z$ and
$(\psi_u)_*v=v$. Therefore $\psi_u (\, c_{u_0}(z,v)(t)\, )$, $t\geq 0$ is the 
$\sigma(u)$ geodesic that begins at $z$ with direction $v$. Thus 
$\psi_u (c_{u_0}(z,v)(t))=c_{u}(z,v)(t)$, for $t\geq 0$.
Consequently $f_u (\,  c_{u_0}(z,v)(t)   \, )=\psi_u (\, c_{u_0}(z,v)(t\, ))$, $t\geq 1$.
This proves $f_u|_{Q\sm W}=(\psi_u)|_{Q\sm W}$ because every point in $Q\sm W$ belongs to some $\sigma(u_0)$ 
geodesic $c_{u_0}(z,v)(t)$.
Now, since $\psi_u$ is at bounded distance from the identity (recall that (iii) above holds
for $\psi$) then $f_u( \, c_{u_0}(z,v)\, )$ is at bounded distance from $c_{u_0}(z,v)$,
thus they define the same point in $\p_\infty$. Therefore
$f_u( \, [c_{u_0}(z,v)\, ] \,  )\, =[\,  c_{u_0}(z,v)\, )\, ]$.
Hence $(f_u)|_{\p_\infty}=1_{\p_\infty}$. This proves the Claim.\\

By means of an orientation preserving homeomorphism $[1,\infty]\ra[0,1]$ we can
identify $[1,\infty ]$ with $[0,1]$. It follows from Claim 3 that 
we can consider $f_u\in P(N)$.
And we obtain, by Claim 4, a continuous map $f:\D^{k+1}\ra P(N)$. 
We choose this identification map to be linear when restricted to the interval
$[r,2r]$ with image the interval $[\frac{1}{3},\frac{2}{3}]$.
The next Claim proves Theorem 1.\\

\noindent {\bf Claim 6.} {\it  $f|_{\bS^k}$ is homotopic to $\iota_N\theta$}.\vspace{.1in}

\noindent {\bf Proof of Claim 6.}  
Let $u\in\bS^k$.
Recall that $\psi_u$ is the identity outside the union of $\bigcup U'_i$ and $U'$
and inside the closed normal geodesic tubular neighborhood of width $r$ of $\alpha_0=\bS^1$
(see (iii) above). In particular $\psi_u$ is the identity on $W$.
>From (iv')  (and (1)) above we have\\

$(\psi_u)|_{U'}=\left[\Phi^Q(\alpha_0, E,r)\theta (u)\right]|_{U'}$, for $u\in\bS^k$. \\

Recall also that each  $U'_i$ is diffeomorphic to $\D^{n-1}\times\R$.
Let $\bar{\alpha}_0$ be the (not necessarily embedded) closed $g$ geodesic which is the image
of $\alpha_0\sbs Q$ by the covering map $Q\ra M$.
Remark that $U_i$ is the $2r$ normal geodesic tubular neighborhood of a lifting
$\beta_i$ of $\alpha\sbs M$ which is diffeomorphic to $\R$. Since $\alpha\sbs M$ is
freely homotopic to the closed geodesic $\bar{\alpha}_0\sbs M$
we have that $\beta_i$ is at finite distance from some embedded
geodesic line which is a lifting of $\bar{\alpha}_0$. Therefore the closure of 
$U_i$ in $Q\cup\p_\infty$ is formed exactly by the two points at infinity
determined by this geodesic line.  Consequently, the closure $\bar{U}_i$ of each $U_i$ 
is homeomorphic to $\D^n$ and intersects $\p_\infty$ in exactly two different points.
Now, applying Alexander's trick to each $\psi|_{\bar{U}_i}$, we obtain an isotopy 
(rel $U'$) that isotopes $\phi_u$ to a map that is the identity outside $U'\sm int\, (W)$, and 
coincides with $\psi_u$ on $U'$, that is, coincides with $\Phi^Q(\alpha_0,E[\frac{1}{3},\frac{2}{3}],r)\theta (u)$
on $U'$. (Note that this isotopy can be defined because the diameters of the closed sets $\bar{U}_i$
in $(Q\sm int\, W)\cup\p_\infty=N\times [1,\infty]$ converge to zero as $i\ra\infty$). Here we refer to any 
metric compatible with the topology of $N\times [1,\infty]$.) Therefore $\psi_u$ is canonically isotopic to a 
map $\vartheta_u$ that is the identity outside $U'$ and on $U'$ coincides with $\Phi^Q(\alpha_0, E,r)\theta (u)$.
In fact $\vartheta_u$ is the identity outside $N\times  [r,2r]\sbs U\sm W\sbs N\times [1,\infty]$. That is,
for $t\in[1,r]\cup [2r,\infty]$, $\vartheta_u( (z,v),t)=((z,v),t)$, $(z,v)\in N$.

On the other hand, we can deform $\theta_u$ to $\theta'_u$, where $\theta'_u$ is the identity on
$N\times (\,  [0,\frac{1}{3}]\cup [\frac{2}{3},1]\, )$ and $\theta'_u((z,v),t)=\theta'_u((z,v),3t-1)$, for
$t\in [\frac{1}{3},\frac{2}{3}]$.
Finally using the identification mentioned before this Claim, we obtain that $\theta'=\vartheta$. This proves 
Claim 6 and Theorem 1.\vspace{.6in}

\noindent {\bf \large  Section 4. Proof of Theorem 2.}\\

First we recall some definitions and introduce some notation. For a compact manifold $M$, the spaces of smooth and topological
pseudo-isotopies of $M$ are denoted by $P^{diff}(M)$ and $P(M)$, respectively. Both $P^{diff}(M)$ and $P(M)$ are groups
with composition as the group operation.
We have stabilization maps $\Sigma :P(M)\ra P(M\times I)$. The direct limit of the sequence 
$P(M)\ra P(M\times I)\ra P(M\times I^2)\ra\dots\,\, $ is called the space of stable topological pseudo-isotopies of $M$, 
and it is denoted by $\cP (M)$. We define $\cP^{diff} (M)$ in a similar way. The inclusion $P^{diff}(M)\ra P(M)$ induces
an inclusion $\cP^{diff}(M)\ra \cP(M)$.
We mention two important facts:
\begin{enumerate}
\item[{1.}]  $\cP^{diff} (-)$, $\cP (-)$ are homotopy functors. 

\item[{2.}] The maps $\pi_k(P^{diff}(M))\ra\pi_k(\cP^{diff} (M))$, $\pi_k(P(M))\ra\pi_k(\cP (M))$
are isomorphisms for $max\{ 2k+9, 3k+7  \}\leq dim\, M$, see \cite{I}.\\
\end{enumerate}

\noindent {\bf Lemma 4.1.} {\it For every $k$ and every compact smooth manifold $M$,  
the kernel and the cokernel of $\pi_k(\cP^{diff}(M))\ra\pi_k(\cP(M)))$ are finitely generated.}\vspace{.1in}

\noindent {\bf Proof.} We have a long exact sequence (see \cite{Hat}, p.12):
$...\ra \pi_{k+1}({\cal{P}}_S(M)))\ra\pi_k(\cP^{diff}(M)))\ra\pi_k(\cP(M)))\ra\pi_k({\cal{P_S}}(M))\ra ...$,
where ${\cal{P}}_S(M)=lim_n\Omega^n\cP(S^nM)$. An important fact here is that
$\pi_*({\cal{P}}_S(M))$ is a homology theory with coefficients in $\pi_{*-1}(\cP^{diff}(*))$.
Since these groups are finitely generated (see \cite{Dw}) the Lemma follows.\\

\noindent Lemma 4.1 together with (2.) imply:\\

\noindent {\bf Corollary 4.2.} {\it For every $k$ and smooth manifold $M^n$  
the kernel and the cokernel of $\pi_k(P^{diff}(M))\ra\pi_k(P(M)))$ are finitely generated for 
$max\{ 2k+9, 3k+7  \}\leq dim\, M$.}\vspace{.1in}

Write $\iota':DIFF((\bS^1\times\bS^{n-2})\times I,\p)\ra P^{diff}(\bS^1\times\bS^{n-2})$.
Since $\iota_{\bS^1\times\bS^{n-2}}:DIFF((\bS^1\times\bS^{n-2})\times I,\p)\ra P(\bS^1\times\bS^{n-2})$ factors
through $\iota'$,  Corollary 4.2 implies that to prove Theorem 2 it is enough to prove:\\

\noindent {\bf Theorem 4.3.} {\it Let $p$ be a prime integer ($p\neq 2$) such that \,$6p-5<n$.
Then for $k=2p-4$ we have that
$\pi_k(DIFF(\bS^1\times\bS^{n-2}\times I,\p ))$ contains $(\Z_p)^\infty$ and
$\pi_k(\iota' )$ restricted to $(\Z_p)^\infty$ is one-to-one.
When $p=2$, $n$ needs to be \, $\geq 10$. 
Also, if $n\geq 14$, then $\pi_1(DIFF(\bS^1\times\bS^{n-2}\times I,\p ))$ contains $(\Z_2)^\infty$ and
$\pi_1(\iota' )$ restricted to $(\Z_2)^\infty$ is one-to-one.}\\

We will need a little more structure. There is an involution ``\, $-$\, '' defined on $P^{diff}(M)$
by turning a pseuso-isotopy upside down.
For $M$ closed we can define this involution easily in the following way. Let $f\in P^{diff}(M)$.
Define $\bar{f}=(\, (f_1)^{-1}\times 1_I\, ) \circ \hat{f}$, where $\hat{f}=r\circ f\circ r$,
$r(x,t)=(x,1-t)$ and $(f_1(x),1)=f(x,1)$. This involution homotopy anti-commutes
with the stabilization map $\Sigma$, hence the involution can be extended to $\cP(M)$.
This involution induces an involution $-:\pi_k(\cP(M))\ra\pi_k(\cP(M))$ at the $k$-homotopy
level. We define now a map $\Xi: P^{diff}(M)\ra P^{diff}(M)$ by $\Xi(f)=f\circ \bar{f}$, and extend this map to $\cP^{diff}(M)$.
We have four comments:
\begin{enumerate}
\item[{i.}] For $f\in P^{diff}(M)$, $\Xi(f)|_{M\times\{ 1\}}=1_{M\times\{ 1\}}$. Therefore $\Xi(f)\in DIFF(M\times I, \p)$.
Hence  the map $\Xi:P^{diff}(M)\ra P^{diff}(M)$ factors through $ DIFF(M\times I, \p)$.

\item[{ii.}] Since $P^{diff}(M)$ is a topological group, for $x\in\pi_k(P(M))$ we have that
$\pi_k(\Xi)(x)=x+\bar{x}$. 

\item[{iii.}] The following diagram commutes
$$\begin{array}{ccccc}  P^{diff}(M)&&\ra&&P^{diff}(M)\\
\downarrow&&&&\downarrow\\
\cP^{diff}(M)&&\ra&&\cP^{diff}(M)
\end{array}$$
\noindent where the horizontal lines are both either ``$\, -\,$'' or $\Xi$. Hence we have an analogous
diagram at the homotopy group level. 

\item[{iv.}] We mentioned in (1.) that $\cP^{diff}(-)$ is a homotopy functor. But the conjugation ``$-$''
defined on $\cP^{diff}(M)$ depends on $M$. In any event, we have that $\cP^{diff}(-)$ preserves the 
conjugation ``$-$'' up to multiplication by $\pm 1$.
\end{enumerate}

Note that (i.) above implies that 
$\pi_k(\Xi): \pi_k(P^{diff}(\bS^1\times\bS^{n-2}))\ra \pi_k(P^{diff}(\bS^1\times\bS^{n-2}))$ 
factors through $\pi_k(\, DIFF((\bS^1\times\bS^{n-2})\times I,\p)\, )$.
Therefore, to prove Theorem 4.3 it is enough to prove:\\

\noindent {\bf Proposition 4.4.} {\it For every $k=2p-4$, $p$ prime integer ($p\neq 2$),  $6p-5<n$, we have that 
$\pi_k(P^{diff}(\bS^1\times\bS^{n-2}))$ contains $(\Z_p)^\infty$.
Also $\pi_1(P^{diff}(\bS^1\times\bS^{n-2}))$ contains $(\Z_2)^\infty$, provided $n\geq 14$,
and $\pi_0(P^{diff}(\bS^1\times\bS^{n-2}))$ contains $(\Z_2)^\infty$, provided $\geq 10$.
Moreover, in all cases above,  $\pi_k(\Xi)$ restricted these subgroups is one-to-one.}\\

By (2.) and  (iii.) to prove
Proposition 4.4 it is enough to prove the following stabilized version:\\

\noindent {\bf Proposition 4.5.} {\it For every $k=2p-4$, $p$ prime integer ($p\neq 2$),  $6p-5<n$, we have that 
$\pi_k(\cP^{diff}(\bS^1\times\bS^{n-2}))$ contains $(\Z_p)^\infty$.
Also $\pi_1(\cP^{diff}(\bS^1\times\bS^{n-2}))$ contains $(\Z_2)^\infty$, provided $n\geq 14$,
and $\pi_0(\cP^{diff}(\bS^1\times\bS^{n-2}))$ contains $(\Z_2)^\infty$, provided $\geq 10$.
Moreover, in all cases above,  $\pi_k(\Xi)$ restricted these subgroups is one-to-one.}\\


Since $\bS^1$ is a retract of $\bS^1\times\bS^{n-2}$, (1.) implies that 
$\pi_k ( \cP^{diff}(\bS^1))$  is a direct summand of $\pi_k ( \cP^{diff}(\bS^1\times\bS^{n-2}))$.
Therefore, by (ii.) and (iv.), to prove
Proposition 4.5 it is enough to prove the following version for $\bS^1$:\\

\noindent {\bf Proposition 4.6.} {\it For every $k=2p-4$, $p$ prime integer, we have that 
$\pi_k(\cP^{diff}(\bS^1))$ contains $(\Z_p)^\infty$.
Also $\pi_1(\cP^{diff}(\bS^1))$ contains $(\Z_2)^\infty$.
Moreover, in these cases, the two group endomorphisms $x\mapsto x+\bar{x}$
and  $x\mapsto x-\bar{x}$ are both one-to-one when restricted to these subgroups.}\\


\noindent {\bf Proof.} 
For a finite complex $X$, Waldhausen \cite{Wa} proved that the kernel of the split
epimorphism  $$\zeta_k:  \pi_k(A(X))\ra\pi_{k-2}(\cP^{diff}(X))$$ is finitely generated.
Recall that the conjugation in $\cP^{diff}(X)$ is defined by turning a pseudo-isotopy
upside down. It is also possible to define a conjugation ``$-$'' on $A(X)$ such that $\zeta_k$ preserves
conjugation up to multiplication by $\pm 1$ (see \cite{V}). The induced map at the
$k$-homotopy level will also be denoted by  ``$-$''.

We recall a result proved in \cite{HKVW}. For a space $X$ we have that $\pi_k(A(X\times\bS^1))$
naturally decomposes as a sum of four terms: 
$$\pi_k(A(X\times\bS^1))=\pi_k(A(X))\oplus\pi_{k-1}(A(X))\oplus\pi_k(N_-A(X))\oplus\pi_k(N_+A(X))$$
\noindent and the conjugation leaves invariant the first two terms and interchanges the last two.\\

The following result is crucial to our argument:\\

\noindent {\bf Theorem ($\, p\,$-torsion of $\pi_{2p-2}A(\bS^1)$)}. {\it For every prime $p$ the
subgroup of $\pi_{2p-2}(A(\bS^1))$ consisting of all elements of order $p$ is isomorphic
to $(\Z_p)^\infty$.}\\

\noindent  (We are grateful to Tom Goodwillie for communicating this result to us.)
Also Igusa (\cite{I1}, Part D, Th 2.1) building on work of Waldhausen \cite {Wa} proved the following:\\

\noindent {\bf Addendum.} {\it $\pi_3A(\bS^1)$ contains $(\Z_2)^\infty$.}\\

\noindent {\bf Remark.} The special case of the $p$-torsion Theorem above, when $p=2$, is also due to Igusa
(see \cite{I1}, Th. 8.a.2).\\

\noindent
Now, take $X=*$ in the decomposition formula above. Recall that Dwyer showed that $\pi_k(A(*))$ is finitely generated for all $k$.
Therefore the Theorem above implies that at least one of the summands $\pi_k(N_-A(*))$, $\pi_k(N_+A(*))$
in the above formula contains $(Z_p)^\infty$, for $k=2p-2$ and contains $(\Z_2)^\infty$ when $k=3$ by the Addendum.
Hence $y\mapsto y+\bar{y}$ and  $y\mapsto y-\bar{y}$, $y\in (\Z_p)^\infty$, are both one-to-one.
Since $\zeta_k:  \pi_k(A(X))\ra\pi_{k-2}(\cP^{diff}(X))$
has finitely generated kernel we can assume 
(by passing to a subgroup of finite index) that $y\mapsto \zeta_k(y+\bar{y})$
and $y\mapsto \zeta_k(y-\bar{y})$, $y\in (\Z_p)^\infty$, are also one-to-one. It follows that
$x\mapsto x+\bar{x}$ and  $x\mapsto x-\bar{x}$, $x\in \zeta_k((\Z_p)^\infty)$, are one-to-one.
Finally, the same argument shows that 
$x\mapsto x+\bar{x}$ and  $x\mapsto x-\bar{x}$, $x\in \zeta_3((\Z_2)^\infty)$, are one-to-one.

F.T. Farrell

SUNY, Binghamton, N.Y., 13902, U.S.A.\\

P. Ontaneda

SUNY, Binghamton, N.Y., 13902, U.S.A.


\begin{thebibliography}{99}









\bibitem{BH} M. Bridson and A. Haeflinger, {\em Metric spaces of 
non-positive curvature},
Springer-Verlag (1999).


   

\bibitem{BK} K. Burns and A. Katok, {\em Manifolds with non-positive curvature}, Ergodic
Theory \& Dynam. Sys. {\bf 5} (1985), 307-317.









\bibitem{DI} R.K. Dennis and K. Igusa, {\em Hochschild homology and the second obstruction for pseudo-isotopies}, LNM
{\bf 966}, Springer-Verlag, Berlin, 1982, pp. 7-58.





\bibitem{Dw} W. Dwyer, {\em Twisted homology stability for general linear groups}, Ann. of Math. (2)
{\bf 111} (1980) 239-251.


\bibitem{EE} C. J. Earle and J. Eells, {\em Deformations of Riemannian surfaces}, LNM {\bf 102}, Springer-Verlag,
Berlin (1969) 122-149.





\bibitem{EL2} J. Eells and L. Lemaire, {\em Deformations of metrics and associated harmonic maps}, 
Patodi Memorial Volume, Geometry and Analysis, Tata Institute, Bombay, 1981, 33-45. 


























\bibitem{Good} T.G. Goodwillie, {\em The differential calculus of homotopy functors}, in Proceedings of the International
Congress of Mathematicians (Kyoto, 1990), pp. 621-630. Math. Soc. Japan, Tokyo, 1991.



\bibitem{H} R. Hamilton, {\em The Ricci flow on surfaces}, Contemporary Mathematics {\bf 71} (1988)
237-261.






\bibitem{Hat} A.E. Hatcher, {\em Concordance spaces, higher simple homotopy theory, and
applications}, Proc. Symp. Pure Math. {\bf 32} (1978) 3-21.










\bibitem{HKVW} T. Hutterman, J.R. Klein, W. Vogell, F. Waldhausen and B. Williams, {\em The ``fundamental theorem''
for the algebraic K-theory of spaces II- the canonical involution}, Journal of Pure and Applied Algebra 
{\bf 167} (2002) 53-82.


\bibitem{I1} K. Igusa, {\em What happens to Hatcher and Wagoner's formula for $\pi_0(\cC (M))$
when to first Postnikov invariant of $M$ is non-trivial?}, LNM {\bf vo. 1046}, Springer-Verlag, Berlin, 1984, pp. 104-177.

\bibitem{I} K. Igusa, {\em Stability Theorems for pseudo-isotopies}, 
K-theory {\bf 2} (1988) 1-355.





















































\bibitem{V} W. Vogell, {\em The involution in algebraic K-theory of spaces}, Algebraic and Geometric Topology,
LNM. {\bf vol. 1126}, Springer, Berlin, 1985, 277-317.

\bibitem{Wa} F. Waldhausen, {\em Algebraic K-Theory of topological spaces I},
Proc. Sympos. Pure Math., {\bf 32} (1978)  35-60. 










\end{thebibliography}
\end{document}